\documentclass[12pt]{article}
\usepackage[dvips]{graphicx}
 
\parskip=\medskipamount

\def\be{\begin{equation}}
\def\ee{\end{equation}}
\def\bea{\begin{eqnarray}}
\def\eea{\end{eqnarray}}
\def\bean{\begin{eqnarray*}}
\def\eean{\end{eqnarray*}}
\def\r#1{(\ref{#1})}
\def\la#1{\label{#1}}
\def\c#1{\cite{#1}}

\title{Commuting Extensions and Cubature Formulae}
\author{Ilan Degani\\
        Departments of Mathematics and Chemical Physics\\
        Weizmann Institute of Science, Rehovot 76100, Israel\\
        {\small email: {\tt ilan.degani@weizmann.ac.il}}\\
        Jeremy Schiff \\
        Department of Mathematics\\ 
         Bar-Ilan University, Ramat Gan 52900, Israel\\
        {\small email: {\tt schiff@math.biu.ac.il}}\\
        David Tannor\\
        Department of Chemical Physics\\
        Weizmann Institute of Science, Rehovot 76100, Israel\\
        {\small email: {\tt david.tannor@weizmann.ac.il}}}
\date{July 2004}

\begin{document}

\maketitle

\begin{abstract}
Based on a novel point of view on $1$-dimensional Gaussian quadrature, 
we
present a new approach to the computation of $d$-dimensional cubature
formulae. It is well known that the nodes of $1$-dimensional Gaussian
quadrature can be computed as eigenvalues of the so-called Jacobi
matrix. The $d$-dimensional analog is that cubature nodes can be obtained
from the eigenvalues of certain mutually commuting matrices. These 
are obtained by extending (adding rows and columns to) certain noncommuting
matrices $A_1,\ldots,A_d$, related to the coordinate operators $x_1,\ldots,x_d$,
in ${\bf R}^d$. We prove a correspondence between cubature formulae and ``commuting
extensions'' of $A_1,\ldots,A_d$, satisfying a compatibility condition which,
in appropriate coordinates,
constrains certain blocks in the extended matrices to be zero.
Thus, the
problem of finding cubature formulae can be transformed to the problem 
of
computing (and then simultaneously diagonalizing) commuting extensions.
We give a general discussion of existence and of the expected size of
commuting extensions and describe our attempts at computing them, as well 
as examples of cubature formulae obtained using the new approach.
\end{abstract}

\topmargin -20mm
\textheight 240mm
\vfill\eject

\section{Introduction}

One of the most elegant topics in numerical analysis is the theory of 
Gaussian quadrature \c{1}. Unfortunately this theory is limited to 
one dimension, and although something is known about generalizations to 
multiple dimensions (see \c{2} for a survey article and many references), 
at the moment there are many more questions than 
answers. The aim of this paper is to 
present a new approach to cubature rules (``cubature'' seems to be the 
name given to the generalization of quadrature to arbitrary dimension).
In classical, one-dimensional Gaussian quadrature, the most widely used
method for computing nodes and weights, developed about 35 years
ago \c{3}, involves solving the eigenproblem 
for a certain tridiagonal matrix (see \c{4} for a recent ``basis independent''
discussion of this). As far as we are aware, no extension
of this to higher dimensions has previously been obtained. Our 
proposed method for 
computing $d$-dimensional cubature formulae involves the construction of 
$d$ matrices (with tridiagonal block structure in suitable bases),
extending these, in a manner we will explain below,
to a set of commuting matrices, and then solving the 
simultaneous eigenproblem for these commuting matrices. 

The main novel step in this process, which follows very naturally from a 
new approach we present to one-dimensional Gaussian quadrature, 
is the need to construct {\em commuting extensions} of a set of 
matrices. We say the $N\times N$ 
matrices $\tilde{A}_1,\tilde{A}_2,\ldots,
\tilde{A}_d$ are $N\times N$ 
commuting extensions of the $n\times n$ matrices 
$A_1,A_2,\ldots,A_d$ (here $N\ge n$) if 
the top left $n\times n$ block in $\tilde{A}_i$ is $A_i$, for 
each $i=1,2,\ldots,d$, and  the matrices 
$\tilde{A}_1,\tilde{A}_2,\ldots,\tilde{A}_d$ pairwise commute.
The idea of commuting extensions is 
very natural, but we do not find any such notion in the 
linear algebra or numerical linear algebra literature 
(see, however, some very similar ideas in a recent, independent, 
poster of Cargo and Littlejohn \c{cl}).
Since we hope the idea of commuting extensions
will find other applications, the first
few sections in this paper explore this
subject without reference to cubature rules.
Section 2 covers basic theory, section 3
discusses a couple of simple algorithms for 
computing commuting extensions, with which we currently 
have very limited success, and section 4
discusses commuting extensions when the matrices $A_i$ take
a special form,  relevant for the study of cubature
rules. 

In section 5 we turn to the theory of cubature rules. Subsection 
5.1 contains a novel approach to one dimensional Gaussian
quadrature, based upon the properties of a 
certain  operator, its eigenvalues and eigenfunctions. 
This serves as the model for all the subsequent discussion. 
In subsection 5.2 we consider the natural extension of 
this approach to multiple dimensions, and prove the central
results of the paper, giving an equivalence
between  odd degree, positive weight cubature rules
and commuting extensions (satisying a 
compatibility condition that will be explained in the sequel) 
of a certain set of matrices. 
Subsection 5.3 includes some simple consequences of this
relationship. One of the key results in the theory
of cubature rules, a lower bound on the number of nodes needed for an
odd degree cubature rule, originating in the work of 
M\"oller \c{Moller}, 
follows from a general result
in the theory of commuting extensions (theorem 2 in section 3). 
Similarly a simple result on the spectra of commuting extensions
(theorem 6 in section 3) gives interesting constraints on the nodes 
in positive weight cubature rules, which we believe have hitherto been
overlooked even in one dimension. 

In section 6 we turn to actual application of 
the commuting extension approach for computing nodes and weights. 
As in section 3, our achievements here are rather limited,
nevertheless they validate our approach, and even give a
few new cubature formulae. 
Section 7 contains a list of open questions.

We close this introduction by mentioning that one of the oldest continuing 
applications of Gaussian quadrature is in quantum mechanics, where it is used, 
in the so-called ``DVR method'' for the computation of matrix elements of 
non-exactly solvable Hamiltonians (see \c{5} for early references and \c{6} 
for  recent reviews). This paper was born out of an attempt to extend the DVR 
method to higher dimensions; other, recent progress on this subject has 
been made by Dawes and Carrington \c{7}. Another interesting perspective 
on DVR can be found in the paper \c{lccmp}. 

\section{Commuting Extensions}

In this section we present the basic theory of commuting extensions.

\noindent{\bf Definition}. We say the $N\times N$ 
matrices $\tilde{A}_1,\tilde{A}_2,\ldots,
\tilde{A}_d$ are  {\em $N\times N$
commuting extensions} of the $n\times n$ matrices 
$A_1,A_2,\ldots,A_d$ (here $N\ge n$) if 
the top left $n\times n$ block in $\tilde{A}_i$ is $A_i$, for 
each $i=1,2,\ldots,d$, and  the matrices 
$\tilde{A}_1,\tilde{A}_2,\ldots,\tilde{A}_d$ pairwise commute.

\noindent{\bf Theorem 1.} Any set of 
matrices admits commuting extensions. 

\noindent{\bf Proof:} We construct explicit commuting extensions
of the $n\times n$ matrices $A_1,A_2,\ldots,A_d$. Take 
\be 
\tilde{A}_1 = \left(
\begin{array}{ccccc}
A_1 & A_2 & A_3 & \ldots & A_d  \\
A_d & A_1 & A_2 & \ldots & A_{d-1}  \\
A_{d-1} & A_d & A_1 & \ldots & A_{d-2}  \\
\vdots & \vdots & \vdots & & \vdots \\
A_2 & A_3 & A_4 & \ldots & A_1 
\end{array}
\right)\ ,\quad
\tilde{A}_2 = \left(
\begin{array}{ccccc}
A_2 & A_3 & A_4 & \ldots & A_1  \\
A_1 & A_2 & A_3 & \ldots & A_d  \\
A_d & A_1 & A_2 & \ldots & A_{d-1}  \\
\vdots & \vdots & \vdots & & \vdots \\
A_3 & A_4 & A_5 & \ldots & A_2 
\end{array}
\right)\ ,\quad
{\rm etc.}
\ee 
More fully, take $\tilde{A}_i$ 
to be a $dn \times dn$ matrix 
which is a  $d\times d$ matrix of $n\times n$ blocks, 
with the $j,k$th block, which we denote $(\tilde{A}_i)_{jk}$,
equal to $A_{i+k-j\ ({\rm mod}~d)}$. Then
\be
(\tilde{A}_i \tilde{A}_{i'})_{jj'} = 
   \sum_{k=1}^d   (\tilde{A}_i)_{jk} (\tilde{A}_{i'})_{kj'}   =  
   \sum_{k=1}^d  A_{i+k-j\ ({\rm mod}~d)} A_{i'+j'-k\ ({\rm mod}~d)}\ .
\ee
Replacing the summation index $k$ by $k'+i'-i\ ({\rm mod}~d)$ 
the second sum becomes 
\be
   \sum_{k'=1}^d  A_{i'+k'-j\ ({\rm mod}~d)} A_{i+j'-k'\ ({\rm mod}~d)}\ 
\ee
which is equal to $(\tilde{A}_{i'} \tilde{A}_{i})_{jj'}$. Thus 
\be
(\tilde{A}_i \tilde{A}_{i'})_{jj'}=
(\tilde{A}_{i'} \tilde{A}_{i})_{jj'}
\ee
for all $j,j'$, i.e. $\tilde{A}_i$ and $\tilde{A}_{i'}$ commute.
$\bullet$

Theorem 1 establishes the existence of commuting extensions, but it is
natural to ask what is the smallest possible dimension  for 
commuting extensions 
of a given set of matrices. To this end we have the following result:

\noindent{\bf Theorem 2.} If $N\times N$ commuting extensions
of the $n\times n$  matrices $A_1,A_2,\ldots,A_d$ exist, then
\be N \ge n + \frac12 {\rm max}_{i,j} {\rm rank}([A_i,A_j])\ . \ee 

\noindent{\bf Proof:} 
Suppose  $\tilde{A}_1,\tilde{A}_2,\ldots,
\tilde{A}_d$ are  $N\times N$  
commuting extensions of the $n\times n$ matrices 
$A_1,A_2,\ldots,A_d$. Write 
\be \tilde{A}_i = 
  \left( 
\begin{array}{cc}
A_i & a_i \\ b_i & \alpha_i 
\end{array}
\right)\ ,  \ee
where the matrices $a_i,b_i,\alpha_i$ have sizes
$n\times(N-n)$, $(N-n)\times n$, $(N-n)\times(N-n)$ respectively.
The top left $n\times n$ block of the equation 
$[\tilde{A}_i,\tilde{A}_j]=0$ gives the requirement
\be [A_i,A_j] + a_i b_j - a_j b_i = 0 \ .\la{1e}\ee
Since the matrices $a_i$ and $b_i$ do not have rank exceeding
$(N-n)$, neither do products of the form $a_i b_j$, and the matrices
$a_i b_j - a_j b_i$ can have rank at most $2(N-n)$. Thus 
\r{1e} can hold only if for each $i,j$ we have
\be {\rm rank}([A_i,A_j]) \le 2(N-n)\ , \ee
and the theorem follows directly. $\bullet$

Unfortunately there is a large gap between the lower bound on $N$
from theorem 2 and the $N$ in the existence proof of theorem 1.  
In practice, it seems that the lower bound of theorem 2 is 
rarely attained, and the $N$ of theorem 1 is much too big.
As we shall see in Section 5, theorem 2 gives rise
to a well known lower bound on the number of points needed for 
a cubature formula, and in that context also the bound can rarely be
attained. 

In addition to not knowing, in general, any way to rigorously
predict the lowest dimension for commuting extensions of a given set
of matrices, we also currently have no way of determining
how many distinct families of commuting extensions of a given dimension exist. By a family
we mean a set of commuting extensions related by conjugation as described
in the following obvious result:

\noindent{\bf Theorem 3.}
If the 
matrices $\tilde{A}_1,\tilde{A}_2,\ldots,
\tilde{A}_d$ are  $N\times N$ 
commuting extensions of the $n\times n$ matrices 
$A_1,A_2,\ldots,A_d$, then so are 
the 
matrices $\tilde{U}\tilde{A}_1 \tilde{U}^{-1}, \tilde{U}\tilde{A}_2 \tilde{U}^{-1},\ldots,
\tilde{U}\tilde{A}_d \tilde{U}^{-1}$, where $\tilde{U}$ is any  matrix of the form
\be \tilde{U} = \left(
\begin{array}{cc}
I_{n \times n} &  0_{n \times (N-n)} \\
0_{(N-n) \times n} & U 
\end{array}
\right)
\ee
with $U$ an invertible $(N-n)\times(N-n)$ matrix. 

To proceed further, and at least get some idea of the size needed
for commuting extensions, we have to resort to parameter counting.
From here on we restrict to the case where the matrices 
$A_i$ and $\tilde{A}_i$ are symmetric, i.e. the case of 
{\em symmetric} commuting extensions of a set of {\em symmetric} 
matrices. Note that except when $d=2$ 
the existence construction of theorem 1 does not guarantee
symmetric commuting extensions. Neither is it clear 
that the lowest dimension commuting extensions of a set of symmetric 
matrices need necessarilly be symmetric. But because the case of 
symmetric commuting extensions of symmetric matrices is relevant
for cubature rules, we restrict our attention to this.

If the matrices $\tilde{A}_i$ are symmetric then we can write
\be \tilde{A}_i = 
  \left( 
\begin{array}{cc}
A_i & a_i \\ a_i^T & \alpha_i 
\end{array}
\right)\ ,  \ee
 where $a_i$ is $n\times(N-n)$ 
and $\alpha_i$ is $(N-n)\times(N-n)$ and symmetric. Thus the number
of free parameters we have in choosing the extensions of the $A_i$ is 
\be
 d\left( n(N-n) + \frac12(N-n)(N-n+1) \right) 
=
\frac12 d(N-n)(N+n+1)\ .
\la{e2}\ee
Let us assume that at least one of the $\tilde{A}_i$, 
say $\tilde{A}_1$, 
has distinct eigenvalues. 
Then all matrices that 
commute with $\tilde{A}_1$ also commute amongst themselves, 
and we just need to check that 
$[\tilde{A}_1,\tilde{A}_i]=0$ for $i=2,\ldots,d$. Since the 
commutator of symmetric matrices is automatically antisymmetric, 
we have
\be \frac12 N(N-1)(d-1) \la{e3}\ee
equations to satisfy. We cannot, however, directly compare the 
number of parameters from \r{e2} with the number of equations 
from \r{e3}, as from theorem 3 we learn that (except when $N=n+1$) 
commuting extensions exist in families. For symmetric commuting 
extensions the matrices $U$ (and thus $\tilde{U}$) in theorem 3
are restricted to be orthogonal. So symmetric commuting extensions
occur in families with $\frac12(N-n)(N-n-1)$ parameters, and the number 
of parameters in choosing extensions should exceed the number of
equations from \r{e3} by at least this amount. Thus we need
\be  \frac12 d(N-n)(N+n+1) \ge 
     \frac12 N(N-1)(d-1) + \frac12(N-n)(N-n-1)\ . \ee
A little rearranging of this inequality gives the condition 
\be  N-n \ge \frac{n(n-1)(d-1)}{2(n+d)} = 
      \frac{d-1}{2}n - \frac{d^2-1}{2} + \frac{d(d^2-1)}{2n}
   + o\left( \frac1{n} \right) \ .
\la{eb1}\ee
If $N$ satisfies this condition we expect to find $N \times N$ 
symmetric commuting extensions.
For comparison, in the explicit commuting extensions 
of theorem 1 (which, however, were not symmetric) we had $N-n=(d-1)n$;
we can clearly expect to do much better than this. 

On occasions it seems that 
parameter counting can be misleading. As an example consider the case of $d=2$.
From theorem 2 we have $N-n \ge \frac12{\rm rank}([A_1,A_2])$
(note that since the commutator of 2 symmetric matrices is antisymmetric, its
rank is always even). The parameter counting argument tells us that we 
should expect 
\be  N-n \ge \frac{n(n-1)}{2(n+2)}\ . \ee
Assuming $[A_1,A_2]$ of maximal rank we have 
\be 
\frac12{\rm rank}([A_1,A_2])=
\left\{ 
\begin{array}{cc}
\frac{n}{2}  &  n {\rm ~even} \\ 
\frac{n-1}{2} & n{\rm ~ odd} 
\end{array}\right. \ .\ee
So by theorem 2
\be 
N-n \ge
\left\{ 
\begin{array}{cc}
\frac{n}{2}  &  n {\rm ~even} \\ 
\frac{n-1}{2} & n{\rm ~ odd} 
\end{array}\right. \ .\ee
We see that when $d=2$ and $[A_1,A_2]$ is of maximal rank 
the inequality from parameter counting is actually 
weaker  than the rigorous one from theorem 2. 
Thus in this case the parameter counting
argument is certainly flawed. But this seems to 
be rather exceptional;
in general it appears that
when it is consistent with the lower bound
of theorem 2, parameter counting gives a better idea of the size
we should expect for commuting extensions. In particular we give the
following example where the lower bound of theorem 2 cannot be attained,
at least with symmetric extensions:

\noindent{\bf Theorem 4.} For $n>5$ 
there exist symmetric $n \times n$ matrices
$A_1,A_2$  with ${\rm rank}\left([A_1,A_2]\right)=2$ and no 
$(n+1)\times(n+1)$ symmetric commuting extensions. 

\noindent The proof proceeds through the following lemma which 
will also be useful later:

\noindent{\bf Lemma 1.} Let $A_1,A_2$ be a pair of symmetric $n \times n$
matrices with ${\rm rank}\left([A_1,A_2]\right)=2$. Let $\{v,w\}$ be a 
basis of ${\rm Im}\left([A_1,A_2]\right)$. If 
$A_1,A_2$ have  $(n+1)\times(n+1)$ symmetric commuting extensions then
the vectors $v,w,A_1v,A_1w,A_2v,A_2w$ are linearly dependent.

\noindent{\bf Proof} (of lemma 1). $(n+1)\times (n+1)$ 
symmetric commuting extensions of $A_1,A_2$ must take the form 
\be \tilde{A}_1=\left( \begin{array}{cc}A_1 & a \\ a^T & \alpha \end{array}
      \right)\ , \qquad
\tilde{A}_2=\left( \begin{array}{cc}A_2 & b \\ b^T & \beta \end{array}
      \right)\ , \la{fsce}\ee
where $a,b$ are $n$-dimensional column vectors and $\alpha,\beta$ are 
scalars. The requirement $[\tilde{A}_1,\tilde{A}_2]=0$ translates 
into the equations 
\bea
[A_1,A_2]+ab^T-ba^T &=& 0\ , \la{re1}\\
A_1b + \beta a - A_2a - \alpha b &=& 0 \la{re2}\ . 
\eea
Finding $(n+1)\times (n+1)$ 
symmetric commuting extensions of $A_1,A_2$ is 
equivalent to finding vectors $a,b$ and scalars $\alpha,\beta$ 
satisfying \r{re1}-\r{re2}. 

Since ${\rm rank}\left([A_1,A_2]\right)=2$, the vectors $a,b$ cannot be 
linearly dependent, for if they were we would have $ab^T-ba^T=0$, giving 
a contradiction with \r{re1}. Thus we can find a vector orthogonal to 
$a$ but not to $b$. 
Applying \r{re1} to this we deduce that $a$ is in ${\rm Im}([A_1,A_2])$. 
Likewise, applying \r{re1} to a vector orthogonal to $b$ but not to $a$,
we see that $b$ is in ${\rm Im}([A_1,A_2])$. Choose 
$\{v,w\}$ to be a basis of ${\rm Im}\left([A_1,A_2]\right)$. 
From the previous remarks it follows that we can write 
\bea   a &=& \lambda v + \mu w  \\
       b &=& \nu v + \rho w \eea
where  $\lambda, \mu, \nu, \rho$ are constants with
$\lambda \rho - \mu \nu \not=0 $. 
Substituting in \r{re2} we have
\be
 (\beta \lambda - \alpha \nu) v 
+ (\beta \mu  - \alpha  \rho) w 
+ \nu A_1 v 
+ \rho A_1 w 
- \lambda A_2 v  
- \mu A_2 w 
= 0 \ .
\ee
Since $\lambda \rho - \mu \nu \not=0 $ we see at once that the 
6 vectors $v,w,A_1v,A_1w,A_2v,A_2w$ must
be linearly dependent. $\bullet$

\noindent{\bf Proof} (of Theorem 4). Let 
\be
(A_1)_{a b} = \left\{ 
\begin{array}{cc}
\lambda_a  & a=b \\
0  &  a\not=b
\end{array}
\right. \ , \qquad
(A_2)_{ab} = \left\{ 
\begin{array}{cc}
\mu_a  & a=b \\
\frac{w_av_b-w_bv_a}{\lambda_a - \lambda_b} & a\not=b 
\end{array}
\right. \ ,
\ee
where the $\lambda_a$ are distinct, the $\mu_a$ are arbitrary, 
and the $v_a,w_a$ are entries of 
two arbitrary linearly independent $n$-dimensional vectors. 
Then  $[A_1,A_2]=wv^T-vw^T$, 
so ${\rm rank}\left([A_1,A_2]\right)=2$, 
and $\{v,w\}$ is a basis of ${\rm Im}\left([A_1,A_2]\right)$. 
By the lemma there 
can only exist $(n+1)\times(n+1)$ 
symmetric commuting extensions of $A_1,A_2$ if 
the 6 
vectors $v,w,A_1v,A_1w,A_2v,A_2w$ are linearly dependent. There is no evident 
reason
why they should be linearly dependent if $n>5$, but to show 
concretely that they 
usually are not, we considered the case $n=6$ where the $\lambda_a$ 
take the 
values $1,2,3,4,5,6$, and then constructed, using Maple, the 
determinant of the 
$6\times 6$ matrix with columns $v,w,A_1v,A_1w,A_2v,A_2w$. The values of the 
$v_a,w_a,\mu_a$
were not fixed. The resulting determinant, which is too long to 
reproduce here, is 
simply a polynomial expression in the eighteen variables $v_a,w_a,\mu_a$, 
and is not indentically zero. 
The case $n>6$ follows trivially from the case $n=6$.  $\bullet$

\noindent
{\bf Note on index conventions:} In discussion of commuting extensions we
start with $d$ matrices   of size $n\times n$ which we 
extend to size $N\times N$. 
For clarity, in most of this paper we adhere to the following index conventions:
\newline Indices $i,j,k$ etc. run from $1$ to $d$.
\newline Indices $a,b,c$ etc. run from $1$ to $n$
\newline Indices $\alpha,\beta,\gamma$ etc. run from $1$ to $N$.

The results up to here all concern the existence and size
of commuting extensions. For the purpose of finding commuting extensions,
to be discussed in the next section, we will use the following:

\noindent{\bf Theorem 5}. The $n\times n$ symmetric matrices 
$A_1,A_2,\ldots,A_d$ admit $N\times N$ symmetric commuting 
extensions if and only if there exist $N \times N$ diagonal matrices 
$\Lambda_1,\Lambda_2,\ldots,\Lambda_d$ and an $n \times N$ matrix $Q$
with orthonormal rows such that 
\be A_i = Q \Lambda_i Q^T \ .\ee

\noindent{\bf Proof}. {\em From the extensions to $\Lambda_i,Q$:} 
If we can find $N\times N$ symmetric commuting 
extensions $\tilde{A}_1,\tilde{A}_2,\ldots,\tilde{A}_d$, then
we can find diagonal matrices 
$\Lambda_1,\Lambda_2,\ldots,\Lambda_d$ and an $N \times N$ 
orthogonal matrix $\tilde{Q}$ such that 
\be \tilde{A}_i = \tilde{Q} \Lambda_i \tilde{Q}^T \ .\ee
The matrix $Q$ in the theorem is comprised of just the first 
$n$ rows of $\tilde{Q}$. 

\noindent{\em From $\Lambda_i,Q$ to the extensions:} A matrix $Q$ as 
described in the theorem can always be extended, by the addition 
of $N-n$ orthonormal rows, to an $N\times N$ orthogonal matrix $\tilde{Q}$. 
(In fact this can be done in many ways, corresponding to the 
freedom described in Theorem 3.) Once such  a $\tilde{Q}$ has 
been constructed the matrices 
$ \tilde{A}_i = \tilde{Q} \Lambda_i \tilde{Q}^T $ are
$N\times N$ symmetric commuting extensions of the $A_i$. 
$\bullet$

\noindent{\bf Note:} The matrix $Q$ in theorem 5 satisfies 
$QQ^T=I_{n\times n}$.  

It is of interest to understand how the spectra of commuting 
extensions (i.e. the entries of the matrices $\Lambda_i$ in Theorem 5)
are related to the spectra of the original matrices $A_i$. The 
following is a first result in this direction. It is  a 
simple consequence of the Sturmian separation  theorem \c{sst}, but 
we offer a direct proof too.

\noindent{\bf Theorem 6}. 
Let $\tilde{A}$ be an $N\times N$ symmetric extension of 
the $n\times n$ matrix $A$. Then the 
smallest eigenvalue of $\tilde{A}$ is less than or equal
to the smallest eigenvalue of ${A}$, and
the largest eigenvalue of $\tilde{A}$ is greater than or equal
to the largest eigenvalue of ${A}$.

\noindent{\bf Proof}.
Diagonalizing $\tilde{A}$ and $A$ we have
\be
\tilde{A} = \tilde{Q} \Lambda \tilde{Q}^T \ , \qquad
A = U D U^T \ ,
\ee
where $Q$ and $U$ are $N\times N$ and $n\times n$ orthogonal matrices,
respectively, and $\Lambda$ and $D$ are
$N\times N$ and $n\times n$ diagonal matrices containing
the eigenvalues of $\tilde{A}$ and $A$, respectively.
Restricting the first equation to its upper $n\times n$ block 
we have also
\be A = Q \Lambda Q^T\ ,\ee
where $Q$ is an $n\times N$ matrix composed of the first 
$n$ rows of $\tilde{Q}$, which are orthonormal. Comparing the 
two formulae for $A$ and using the orthogonality of $U$ we 
have 
\be D = (U^TQ) \Lambda (U^TQ)^T\ . \la{evrel}\ee
The matrix $U^TQ$, like $Q$, has orthonormal rows; in 
particular we have
\be \sum_{\alpha=1}^N (U^TQ)_{a\alpha}^2
  = 1 \ , \qquad  a=1,\ldots,n\ .
\ee
The diagonal entries of \r{evrel} read
\be D_a = \sum_{\alpha=1}^N (U^TQ)_{a\alpha}^2 \Lambda_{\alpha}
    \ , \qquad  a=1,\ldots,n\ .
\ee
Thus each eigenvalue $D_a$ of $A$ is a  
generalized average of the eigenvalues
$\{\Lambda_\alpha\}$ of $\tilde{A}$ with respect to a set of 
positive weights $\{(U^TQ)_{a\alpha}^2\}$. It follows at once that 
it is not possible that any of the $D_a$ be either smaller than, or
greater than, all of the $\Lambda_\alpha$. $\bullet$

\section{Computing Symmetric Commuting Extensions}

The most obvious approach to computing commuting extensions
is simply to treat the unknown entries
in the extended matrices $\tilde{A}_i$ as variables, and to consider
the conditions $[\tilde{A}_i, \tilde{A}_j]=0$ as equations in these
variables. In the generic case (generically we should
expect the $\tilde{A}_i$ 
have distinct eigenvalues) it will be sufficient to look 
at the equations just for one particular value of $i$. If $N - n > 1$, 
then by theorem 3 we expect continuous
families of extensions, this freedom can be exploited to fix some of
the variables. The system of equations we obtain will be quadratic in the 
unknown variables, and can be tackled by standard methods for solving 
systems of equations. We have done some initial experiments with this
approach in the case $d=2$, 
attempting to solve the system of quadratic equations 
1) by integrating the gradient flow $v'=-\nabla || [ \tilde{A}_1(v), 
\tilde{A}_2(v)]||^2$ (here $v$ denotes the variables added to form
the commuting extensions), and 2) using Newton's method. 
The results are very varied; for some pairs of moderate-sized matrices 
there is reasonable convergence, but in other cases there are signs
of extreme ill-coniditioning (very low gradients in the case of gradient
flow, almost singular Jacobian in Newton's method). Some of the 
cubature related results we will present in section 6 were obtained by the
integration of the gradient flow with Frobenius norm of the commutator;
some technical details of these calculations can be found in 
\c{ilanthesis}.

In this section we focus on a different approach 
to computing commuting extensions,
based on theorem 5 from section 2, and thus
restricted to the case of symmetric extensions for symmetric matrices. 
We attempt to construct the matrices $\Lambda_i$ and $Q$ introduced in the 
theorem by minimization of 
\bea 
S(Q,\Lambda_1,\ldots,\Lambda_d) 
&=& \frac12{\rm Tr} \left( \sum_{i=1}^d (A_i - Q\Lambda_i Q^T)^2 \right)  
   \nonumber\\
&=&  \frac12\sum_{i=1}^d \sum_{a,b=1}^n \left((A_i - Q\Lambda_i Q^T)_{ab}
        \right)^2 \ .
\la{Sdef}\eea
Consider first the variation of this with respect to the entries of
$\Lambda_i$ (which is diagonal). It is straightforward to check that
\be
\frac{\partial S}{\partial (\Lambda_i)_{\alpha\alpha}}
=  (Q^T(A_i - Q\Lambda_i Q^T)Q)_{\alpha \alpha} \qquad
\left\{
\begin{array}{c}
i=1,\ldots,d \\
\alpha= 1,\ldots,N
\end{array}
\right. \ ,
\ee 
implying that the entries of $\Lambda_i$ must satisfy 
\be \sum_{\beta=1}^N (Q^TQ)_{\alpha \beta}^2 (\Lambda_i)_{\beta\beta}
     = (Q^T A_i Q)_{\alpha \alpha} \qquad
\left\{
\begin{array}{c}
i=1,\ldots,d \\
\alpha= 1,\ldots,N
\end{array}
\right. \ .
\la{findL}\ee
We see that if we can construct $Q$, we can, using the last formula, 
also construct $\Lambda_i$, at least assuming invertibility of the 
matrix with entries  $(Q^TQ)_{\alpha \beta}^2$. 

To build $Q$, or more precisely $\tilde{Q}$, the $N\times N$ orthogonal
matrix that is an extension of the $n\times N$ matrix $Q$ by addition
of $(N-n)$ more rows, we use a sequence of Jacobi rotations (not using 
the rotations that just mix the last $(N-n)$ rows). More fully: assuming 
we have an initial guess $\tilde{Q}_0$ for the matrix $\tilde{Q}$, we 
construct a new
guess in the form $\tilde{Q}=R(\theta)\tilde{Q}_0$ where $R(\theta)$ 
is a Jacobi rotation of
2 specified rows of $\tilde{Q}_0$ through angle $\theta$. $\theta$ is chosen 
to minimize $S$, where in $S$ we use equation \r{findL} to determine 
the $\Lambda_i$. 
The minimization is done numerically, there does not seem to be 
an explicit formula for $\theta$ available.

Note that this algorithm will actually solve a more general problem 
than that of finding commuting extensions. Suppose we use a value of 
$N$ for which there are no commuting extensions. We can still compute
a minimum of $S$, which means we will find $Q$ and $\Lambda_i$ such 
that  $A_i\approx Q\Lambda_i Q^T$. Extending $Q$ to an $N\times N$
orthogonal matrix $\tilde{Q}$ we will have commuting matrices
$\tilde{A_i}= \tilde{Q}\Lambda_i \tilde{Q}^T$ which are {\em approximately} 
extensions of the $A_i$, i.e. we will be computing commuting 
approximate extensions. In the case $N=n$ we will be computing
commuting approximations to the matrices $A_i$ 
without increasing the dimension.
The question  of the existence and computation of 
commuting matrices which approximate a given 
set of matrices, with commutators of small norm, has a substantial
mathematical history, having first been asked, apparently, in \c{a}. 
For some recent references see \c{b}. For the case of
commuting approximations ($N=n$), the numerical approach
we have just outlined was given in \c{c}, though unlike in
our more general case, it seems that in this case there is an 
explicit formula for the choice of rotation angle $\theta$ at each stage. 
This method is also advocated for numerical simultaneous diagonalization of 
commuting matrices. There are applications of the algorithm 
in statistics and  in signal processing \c{d}, and it
is also employed in \c{7}.

It remains to report on some numerical experiments. In all cases we 
implemented
the algorithm starting with various random choices of $\tilde{Q}$
and applying sweeps of Jacobi rotations of all possible pairs of rows. We attempted
to compute symmetric commuting extensions (and approximate
extensions) for 4 specific pairs of symmetric
$6\times 6$ matrices; the findings are
supported by many other numerical experiments as well. 
The specific matrices used can be found on the internet at \newline
\centerline{\tt http://www.math.biu.ac.il/$\sim$schiff/commext.html.} 
\newline
The pairs of matrices chosen covered the following 4 cases:
\begin{enumerate}
\item ${\rm rank}([A_1,A_2])=2$, $7\times7$ symmetric commuting extensions 
known to exist. 
\item ${\rm rank}([A_1,A_2])=2$, $7\times7$ symmetric commuting extensions 
known not to exist (see lemma 1, section 2).
\item ${\rm rank}([A_1,A_2])=4$, $8\times 8$ symmetric commuting extensions 
known to 
exist, this being the lowest dimension allowed by theorem 2, section 2.
(In fact  presumably many commuting extensions exist, as even after 
accounting for the invariance of theorem 3, section 2, there are more 
parameters than equations).
\item ${\rm rank}([A_1,A_2])=6$, $9\times9$ symmetric commuting extensions 
known to 
exist, this being the lowest dimension allowed by theorem 2, section 2.
(Once again presumably many commuting extensions exist.)
\end{enumerate}

In each case we not only ran the algorithm with $N$ large enough
that we expected to find commuting extensions, but also
with  all other $N\ge 6$ smaller than this (thus computing 
commuting approximations for $N=6$ and 
commuting approximate extensions for $N>6$ but too small
for commuting extensions).  Our observations can
be summarized as follows:

\begin{enumerate}
\item For the actual computation of commuting extensions ($N=7$ in case 1,
$N=8$ in cases 2 and 3, $N=9$ in case 4) the positive result is that in each case
the algorithm was observed to converge. After initial stabilization
the value of $S$ was observed to drop off according to 
\be \ln S = a - bk\ , \ee
where $k$ is the sweep number, and $a$ and $b>0$ are (case-dependent) constants.
This behavior was observed over many orders of magnitude. 
The problematic result is that {\em except in case 2} the value of $b$ was found to
be small, to the extent that obtaining even single precision accuracy required 
many thousands of sweeps. This behavior is not particularly 
surprising given that we are in practice doing a multidimensional optimization 
by sequential one-dimensional optimizations, with 
a fixed choice for the search directions. Remarkably in case 2 a large value 
of $b$ was observed (and this persisted for other choices of matrices with
rank 2 commutator but for which we were searching for extensions with $N=8$).
\item For the computation of commuting approximations ($N=n$), in all cases 
the algorithm converged reasonably quickly, as reported in the literature
\c{c}, though there were noticeable differences in the rate of convergence,
with case 1 being substantially slower than all the other cases. 
 No problems with local minima were observed, different initial choices of $\tilde{Q}$
produced the same minimum value of $S$. Testing of this was 
limited, though.
\item In the search for commuting approximate extensions ($N=7$ in cases
2 and 3, $N=7,8$ in case 4) a  variety of behaviors were observed.
For $N=7$ in case 2 the search was, by far, 
the slowest performed, but
it did ultimately converge. In case
4 the $N=8$ search was noticeably slower than the $N=7$ search, and for
$N=8$ we noticed convergence to different values of $S$, indicating the 
presence of local (approximate) minima. In all cases the values of $S$ 
achieved were at most a few orders of magnitude smaller than the values 
obtained when $N=n$. This gives hope that in the general case, when we 
do not know what the smallest dimension for commuting extensions is, 
we might be able to detect it by minimizing $S$ for different $N$, but
a better minimization algorithm than the current one will certainly be 
necessary.  
\end{enumerate}
It is clear from our results that a lot more work is
necessary on the topic of computing commuting extensions. 

\section{A Special Case Of Commuting Extensions}

We turn to consideration of 
a special case of commuting extensions that turns
out to be relevant for cubature formula.
It is characterized by 2 conditions:
First, the symmetric $n\times n$ 
matrices $A_i$ for which we wish to
find commuting extensions have tridiagonal block form 
\be
A_i = \left(
\begin{array}{cccccc}
\alpha_{i1}  & a_{i1}      &  0            &   \ldots   & 0  & 0 \\
a_{i1}^T     & \alpha_{i2} &  a_{i2}       &   \ldots   & 0  & 0 \\
0            & a_{i2}^T    &  \alpha_{i3}  &   \ldots   & 0  & 0 \\
\vdots       &    \vdots   &   \vdots      &            &  \vdots  & \vdots  \\
0            &    0        &   0           &   \ldots   & \alpha_{i(r-1)} & a_{i(r-1)} \\
0            &    0        &   0           &   \ldots   & a_{i(r-1)}^T & \alpha_{ir} 
\end{array}
\right)\ .
\la{trid}
\ee
Here $\alpha_{i1},\alpha_{i2},\ldots,\alpha_{ir}$ are symmetric square matrices 
of sizes $n_1\times n_1,n_2\times n_2,\ldots, n_{r}\times n_{r}$ respectively,
where $n_1+n_2+\ldots+n_r=n$. The matrices $a_{i1},a_{i2},\ldots,a_{i(r-1)}$ are of
size $n_1\times n_2,n_2\times n_3,\ldots, n_{r-1}\times n_{r}$ respectively. 
The second condition we impose is that  the 
commutator matrices $[A_i,A_j]$ all vanish except for a single block in the bottom right
hand corner, of size $n_{r}\times n_{r}$. 

Let us seek symmetric commuting extensions with the matrices $\tilde{A}_i$,
of size $N\times N$, also taking tridiagonal block form, that is 
\be
\tilde{A}_i = \left(
\begin{array}{ccccccc}
\alpha_{i1}  & a_{i1}      &  0            &   \ldots   & 0  & 0  & 0 \\
a_{i1}^T     & \alpha_{i2} &  a_{i2}       &   \ldots   & 0  & 0 & 0 \\
0            & a_{i2}^T    &  \alpha_{i3}  &   \ldots   & 0  & 0  & 0 \\
\vdots       &    \vdots   &   \vdots      &            &  \vdots  & \vdots & \vdots \\
0            &    0        &   0           &   \ldots   & \alpha_{i(r-1)} & a_{i(r-1)} 
                                  & 0 \\
0            &    0        &   0           &   \ldots   & a_{i(r-1)}^T & \alpha_{ir} 
                                  & a_i \\
0            &    0        &   0           &   \ldots   & 0 & a_{i}^T & \alpha_{i}
\end{array}
\right)\ ,
\ee
where the new blocks $\alpha_i$ are of size $(N-n)\times(N-n)$ and 
are symmetric, and the $a_i$ are of size $n_r\times (N-n)$. 

The questions we wish to ask are (1) what are the equations that the new blocks 
$\alpha_i,a_i$ have to satisfy? and (2) how large  need $N$ be for us to 
have a hope that such extensions exist? 
As in section 2, we assume that
$\tilde{A}_1$ has distinct 
eigenvalues, so we need only check that $\tilde{A}_1$ commutes 
with the $d-1$ matrices $\tilde{A}_2,\ldots,\tilde{A}_d$ and this guarantees
that all the $\tilde{A}_i$ mutually commute.  A brief calculation, using the fact
that the commutators $[A_1,A_i]$ are zero except for a single block, gives 
the following conditions:
\be
\left.
\begin{array}{rcl}
a_{1(r-1)}a_i - a_{i(r-1)}a_1  &=&  0 \\
a_{1(r-1)}^Ta_{i(r-1)}-a_{i(r-1)}^Ta_{1(r-1)} 
  +\alpha_{1r}\alpha_{ir}-\alpha_{ir}\alpha_{1r}
  +a_1a_i^T - a_ia_1^T &=& 0 \\
\alpha_{1r}a_i -\alpha_{ir}a_1 + a_1\alpha_i-a_i\alpha_1 &=&0 \\
a_1^Ta_i- a_i^Ta_1 + \alpha_1 \alpha_i - \alpha_i \alpha_1 &=& 0
\end{array}
\right\} 
\qquad i=2,\ldots,d\ .
\la{nes} \ee
Of these four equations for each $i$, the first is of size $n_{r-1}\times (N-n)$, 
the second is of size $n_r\times n_r$ and antisymmetric, 
the third is of size $n_{r}\times (N-n)$
and the fourth is of size $(N-n)\times (N-n)$ and antisymmetric. Thus there is a total
of 
\be
(d-1)\left(n_{r-1}(N-n)
+\frac12 n_r(n_r-1)
+ n_{r}(N-n)+ \frac12 (N-n)(N-n-1) \right)
\ee
equations to be satisfied. The number of variables available in the 
$a_i$ and $\alpha_i$ is 
\be d\left(n_r(N-n)+\frac12(N-n)(N-n+1)\right)\ . \ee
The system of equations \r{nes}, has an invariance 
\be
a_i \rightarrow a_i g \ , \qquad
\alpha_i \rightarrow g^T \alpha_i g \ , \qquad
i=1,\ldots,r
\ee
where $g$ is an $(N-n) \times (N-n)$ orthogonal matrix. Thus it is not sufficient
that the number of variables simply exceed the number of equations to
be solved, it must exceed the number of equations to be solved by at 
least $\frac12(N-n)(N-n-1)$ to give a full family of solutions. Thus 
we can expect solutions provided:
\bea
& d\left(n_r(N-n)+\frac12(N-n)(N-n+1)\right) \ge 
    \frac12(N-n)(N-n-1) ~~~~~~~~~~~~~~~~~~~~~~~~~~~~&  \nonumber \\
& + (d-1)\left(n_{r-1}(N-n)
+\frac12 n_r(n_r-1)
+ n_{r}(N-n)+ \frac12 (N-n)(N-n-1) \right) &
\eea
Simplifying this gives
\be
N-n \ge \frac{n_r(n_r-1)}{2\left( \frac{n_r+d}{d-1} - n_{r-1} \right)} \ ,
\la{nbd}\ee
where we have made the assumption that the denominator on the right hand side
is positive. 

Note the right hand side in \r{nbd} does not depend on the total dimension, $n$,
of the matrices $A_i$, but just on 
$n_{r-1}$ and $n_r$.  Thus the size of the extensions needed in this special
case may well be substantially smaller than in general. The application 
to cubature of this special form of commuting extensions  will become clear in section 
5.3. 

\section{Cubature Formulae and Commuting Extensions}
Much of the contents of this section, with some additions, are also discussed in 
\cite{ilanthesis}.
In section 5.1 we give a novel presentation on the subject of Gaussian quadrature;
in section 5.2 this is extended to the case of multidimensional cubature. Section
5.3 discusses some practical aspects and consequences of the results of 5.2.

\subsection{A Novel Approach to Gaussian Quadrature}
In Gaussian quadrature we wish to find $q+1$ nodes $x_0,\ldots, x_q$ and 
$q+1$ weights $w_0,\ldots,w_q$ such that the quadrature rule
\be \int_\Omega  w(x) f(x) dx \approx \sum_{i=0}^q  w_i f(x_i) \ee
is exact whenever $f(x)$ is a polynomial of degree at most $2q+1$. 
Here $\Omega$ is some interval or union of intervals and $w(x)\ge 0$ is a suitable 
weight function. {\em Throughout this paper we only consider quadrature 
and cubature rules with positive weights, i.e. $w_i>0$.} 
 
\noindent
{\bf Note:} In this subsection there is no mention of commuting extensions so
we allow ourselves to break our index conventions. In this subsection alone 
indices $i,j,k$ run from $0$ to $q$. 

Denote by ${\cal P}_q$ the space of polynomials of degree at most $q$ with 
the inner product 
\be 
\langle a | b \rangle = \int_\Omega w(x) a(x) b(x) dx \ , \qquad 
   \forall a,b \in {\cal P}_q \ .
\ee
Let $\Pi_q$ be the projector from ${\cal P}_{q+1}$ onto ${\cal P}_q$ parallel to its orthogonal
complement ${\cal P}_q^\perp$ $i.e.$ the obvious orthogonal projection onto ${\cal P}_q$ 
with respect to the inner product above.
We define the operator $\chi: {\cal P}_q \to {\cal P}_q$ by
$\chi p = \Pi_q x p$ for all $p \in {\cal P}_q  $. Since $\chi$ is self adjoint there is an 
orthonormal basis of 
${\cal P}_q$ consisting of eigenfunctions $\{ u_i \}$ of $\chi$, $\chi u_i = \Lambda_i u_i$.
Associated with $\chi$ there is a symmetric bilinear form $X$ on ${\cal P}_q$
defined by $X(a,b) = \langle a | \chi b \rangle$, or equivalently
\be 
 X(a,b) = \int_\Omega w(x) a(x) x b(x) \ d x \ , \qquad     \forall a,b \in {\cal P}_q \ .
\ee
$X$ is diagonalised in the basis $\{ u_i \}$, $X(u_i , u_j) = \Lambda_i \delta_{ij}$.

We prove below that the eigenvalues $\{ \Lambda_i \}$ of $\chi$ provide nodes for a Gaussian
quadrature formula of degree $2 q + 1$. 
Our treatment is the reverse of the classical presentation of Gaussian quadrature 
\cite{1}, \cite{3}, see the explanation after theorem 8.

We need the following remarkable lemma:

\noindent{\bf Lemma 2.}(The $\delta$ lemma) Let $p$ be an arbitrary polynomial in ${\cal P}_{q+1}$. 
Then
\be \langle p|u_i \rangle =  \langle 1|u_i \rangle p(\Lambda_i)   \ , \label{ProjLemm1} \ee
$i.e.$ the  inner product of $p$ with $u_i$ is determined, up to normalization,  
by evaluation of $p$ at $\Lambda_i$.

\noindent{\bf Note:} In this lemma $p$ is allowed to be in ${\cal P}_{q+1}$.

\noindent{\bf Proof.} We prove recursively for $j$ that 
\be \langle x^j|u_i \rangle 
          = \langle 1|u_i\rangle \Lambda_i^j  \ ,  \qquad j=0,\ldots,{q+1}\ . \label{ProjMono} \ee
For $j=0$ the statement is trivial. For
$j>0$, 
\be \langle x^j|u_i \rangle = \langle \chi x^{j-1}|u_i \rangle = 
   \langle  x^{j-1}|\chi u_i \rangle = \Lambda_i  \langle  x^{j-1}| u_i \rangle  \  .       \ee
This provides the recursive step proving (\ref{ProjMono}),
the full result follows by linearity. $\bullet$

\noindent
$\chi$ is an approximation of the operator $x$ and 
this lemma shows that $\{ u_i \}$,
the eigenfunctions of $\chi$, share with $\delta$ functions, which are ``eigenfunctions'' of $x$,
the property that projection of a function on either is done by its evaluation at the appropriate
eigenvalue. This similarity is the reason for the name we give to the $\delta$ lemma.
 
With this lemma it is almost immediate to prove the main result we want:

\noindent{\bf Theorem 7}. Let $f$ be a polynomial of degree at most
$2q+1$. Then 
\be \int_\Omega w(x) f(x) dx =
  \sum_{i=0}^q \langle 1|u_i\rangle^2  f(\Lambda_i) \ , \ee 
i.e. the quadrature rule
\be \int_\Omega  w(x) f(x) dx \approx 
   \sum_{i=0}^q  \langle 1|u_i\rangle ^2   f(\Lambda_i) \ee
is exact of degree $2q+1$. 

\noindent{\bf Proof}. Again we prove the result for $f(x)=x^j$, $j=0,\ldots,2q+1$,
the full result follows by linearity. 
For $j\ge 1$, choose integers 
$n_1,n_2$ between $0$ and $q$ such that $j=n_1+n_2+1$. We then have
\be
\int_\Omega w(x) x^j dx = X( x^{n_1}, x^{n_2} ) 
 =  X \left(  \sum_{k=0}^q \langle x^{n_1} | u_k\rangle u_k ,
             \sum_{i=0}^q \langle x^{n_2} | u_i\rangle u_i
     \right)   
 = \sum_{i=0}^q  \langle 1|u_i\rangle^2 \Lambda_i^j \ .
\ee
In the last step we have used the $\delta$ lemma twice. 

For $j=0$ observe that
\be
\int_\Omega w(x)   dx  = 
\langle 1|1\rangle = 
\left\langle  \sum_{k=0}^q \langle 1 | u_k\rangle u_k \left|
         \sum_{i=0}^q \langle 1 | u_i\rangle u_i 
   \right.\right\rangle
= \sum_{i=0}^q  \langle 1|u_i\rangle^2  \ ,\ee
by the orthonormality of the $u_i$. $\bullet$

Theorem 7 relates the spectrum of $\chi$ to the nodes of Gaussian quadrature. 
It is in fact easy to prove other facts in the theory of Gaussian quadrature using our approach. 
For example, to see that none of the weights vanish just take $p = u_i$ in the $\delta$ lemma.
To see that when $\Omega$ is a single interval $[a,b]$ the nodes must be in its interior 
of the interval just observe that
\be
 b - \Lambda_i = \int_a^b w(x) (b - x) u_i(x)^2 \ d x > 0, \qquad
 a - \Lambda_i = \int_a^b w(x) (a - x) u_i(x)^2 \ d x < 0 \ .
\ee
We also obtain the following widely known characterization of the nodes:

\noindent{\bf Theorem 8}. The nodes $\Lambda_i$ are roots of any nontrivial 
degree $q+1$ polynomial orthogonal to ${\cal P}_q$. 

\noindent{\bf Proof}. Let $p$ be a nontrivial 
degree $q+1$ polynomial orthogonal to ${\cal P}_q$. Then using
the $\delta$ lemma, $0=\langle p|u_i \rangle= \langle 1|u_i\rangle
p(\Lambda_i)$. Since $\langle 1|u_i\rangle$ is nonzero, this 
gives $p(\Lambda_i)=0$. $\bullet$ 

Our presentation on Gaussian quadrature is the reverse of that in
\cite{1}, \cite{3}. 
The starting point in \cite{1}, \cite{3} is that the nodes in degree 
$2 q + 1$ Gaussian quadrature are roots of
the degree $q+1$ polynomial $p$ from theorem 8; it is then shown that 
the eigenvalues of a matrix representation of $\chi$ are equal to
these.
Here we have gone in the other direction; without {\em a priori} 
assumption of the existence of a Gaussian quadrature
formula, we have shown that the eigenvalues of $\chi$ 
are quadrature nodes and as a consequence of the $\delta$ lemma we also 
obtain the fact that they are
roots of the degree $q+1$ polynomial $p$ from theorem 8.

As we shall see in the next subsection,
our approach to Gaussian quadrature  
allows a generalization to higher dimensions. 
However, in the generalization of the $\delta$ lemma the 
polynomial $p$ is restricted to ${\cal P}_q$, not ${\cal P}_{q+1}$,
and this means that while 
theorem 7 can be generalized, theorem 8 cannot, at least not immediately. 
So the characterization
of cubature nodes as roots of a polynomial (or a set of polynomials) seems
to be lost. However the characterization of nodes
as eigenvalues persists, as we now set out to show.

\subsection{Generalization to Cubature}

We denote by ${\cal P}_q  $ the $ n = \pmatrix{d+q \cr d\cr}$ dimensional 
vector space
of polynomials in $d$ variables $x_1, \ldots, x_d$ of total degree up to $q$
(total degree is defined by
${\rm degree}(x_1^{m_1}x_2^{m_2}\ldots x_d^{m_d})=m_1+m_2+\ldots+m_d$).

An $N$-point $d$-dimensional cubature formula
\be \int_\Omega  w(x) f(x) d^dx \approx 
   \sum_{\alpha=1}^N  w_\alpha f(x_\alpha) \ee
is said to be of degree $D$ if it 
is exact whenever $f(x)$ is a polynomial of total degree at most $D$ 
and non-exact for at least one polynomial of total degree $D+1$.
Here $\Omega$ is a suitable region in ${\bf R}^d$ and $w(x)\ge 0$ a suitable 
weight function. The weights $w_\alpha$ are assumed positive. 

We supply ${\cal P}_q$ with the inner product
\be 
\langle a | b \rangle = \int_\Omega w(x)  a(x)  b(x) \ d^d x \ , \qquad 
   \forall \ a, \ b \in {\cal P}_q \ ,
\la{ied}\ee
and define $\Pi_q$, the orthogonal projection operator from  ${\cal P}_{q+1}$ onto  ${\cal P}_q$ with respect
to the above inner product,
in the obvious way. We can then define $d$ self adjoint operators $\chi_1, \ldots, \chi_d$ on ${\cal P}_q$
by
\be
 \chi_i  p = \Pi_q x_i p \ , \qquad \forall p \in {\cal P}_q \ ,
\ee
with related symmetric bilinear forms $X_i : {\cal P}_q \times {\cal P}_q \to {\bf R}$,
\be 
 X_i(a,b) = \langle a | \chi_i b  \rangle = \int_\Omega w(x) a(x) x_i b(x) \ d^d x \qquad
 \forall a,b \in {\cal P}_q \ . \label{qfd}
\ee
Generally $[\chi_i, \ \chi_j] \neq 0$ so we can not find a basis of ${\cal P}_q$ in
which all the $\chi_i$ are simultaneously diagonalised, and we do not have a direct analog of the
one-dimensional case in which the eigenvalues of the single operator $\chi$ served as quadrature
nodes. We shall show, however, that there is a
correspondence between cubature rules and spectra of certain commuting extensions of 
matrix representations of the operators $\chi_1,\ldots, \ \chi_d$.
As a first step towards this we prove the following:

\noindent{\bf Theorem 9.} Let the $n \times n$ matrices $A_1,\ldots,A_d$ 
be the representations of the 
operators $\chi_1,\ldots,\chi_d$ in an arbitrary orthonormal
basis $\{e_a\}$ of ${\cal P}_q$ 
(so $(A_i)_{ab}= \langle e_a  |  \chi_i  e_b \rangle $). 
Suppose that for the region $\Omega$ and weight function $w(x)$ we have a 
degree $2q+1$, $N$ point cubature rule with positive weights. Then
there exist $N\times N$ symmetric commuting extensions of $A_1,\ldots,A_d$.

\noindent{\bf Proof.} 
Suppose the cubature rule takes the form
\be \int_\Omega  w(x) f(x) d^dx \approx 
   \sum_{\alpha=1}^N  w_\alpha f(x_\alpha) \ .\ee
Then, since all integrands are of degree at most $2q+1$, we have
\bea
\delta_{ab} &=& 
\langle e_a|e_b \rangle =
\int_\Omega w(x) e_a(x) e_b(x) d^dx= 
\sum_{\alpha=1}^N w_\alpha e_a(x_\alpha)e_b(x_\alpha) \ , \la{m1}\\  
(A_i)_{ab} &=& \langle e_a  |  \chi_i  e_b \rangle =
\int_\Omega w(x) e_a(x) x_i e_b(x) d^dx = 
\sum_{\alpha=1}^N w_\alpha e_a(x_\alpha)(x_\alpha)_i e_b(x_\alpha) \ .
\la{m2}\eea
Define the $n\times N$ matrix $Q$ by
\be
Q_{a\alpha} = \sqrt{w_\alpha} e_a(x_\alpha)\ , 
\ee
and $N\times N$ diagonal matrices $\Lambda_1\,\ldots,\Lambda_d$ with diagonal entries 
$(\Lambda_i)_{\alpha \alpha} = \Lambda_{i \alpha}$,
\be
\Lambda_{i \alpha} =  (x_\alpha)_i \ .  \label{EigXconnect}
\ee
Equations \r{m1}-\r{m2} read 
\be  I_{n \times n}=QQ^T \ , \qquad  
     A_i = Q\Lambda_i Q^T  \ .      \label{QLamRel}
\ee
Using theorem 5 in section 2 we conclude that $A_1,\ldots,A_d$ have 
$N \times N$ symmetric commuting extensions.
$\bullet$

It is natural to ask whether the matrix commuting extensions of 
theorem 9 are representations of
commuting extensions of the operators $\{ \chi_i \}$ in some $N$ dimensional space of 
functions $V$ which
includes
${\cal P}_q$ as a subspace.
In other words, is it possible to extend the basis $\{ e_a  \}$ of ${\cal P}_q$ 
to an orthonormal basis of $V$ by adding $N-n$ orthonormal functions 
$e_{n+1}, \ldots, e_{N}$ 
in such a way that the $N \times N$ matrices $( \tilde{A}_i )_{\alpha, \beta} = \langle e_\alpha |  x_i e_\beta \rangle
= \int_\Omega w(x) e_\alpha(x)   x_i e_\beta(x) \ d^d x$ commute?
Unfortunately, in all but the simplest cases, we could not find, nor prove the existence of, such functions
$e_{n+1}, \ldots, e_{N}$. 
However, even though we can not view the matrix commuting extensions of theorem 9 as representations of operator
extensions of the $\chi_i$, 
they do satisfy a certain compatibility condition with
the $\chi_i$ which we prove in theorem 10.

Let us introduce the $N$ dimensional vector space $V (={\bf R}^N)$, with the standard inner product, whose elements we denote in bold face. ${\cal P}_q$ is mapped to a subspace of $V$ by the inclusion operator $\iota : {\cal P}_q \to V$
which is defined by $\iota e_1 = {\bf e_1}, \ldots, \iota  e_n = {\bf e_n}$, where $\{ {\bf e_a} \}$ are the first
$n$ members of the standard basis of $V$. 
Extend $\{ {\bf e_a} \}$ to an orthonormal basis $\{ {\bf e_\alpha} \}$ of $V$ by adding 
any orthonormal basis $\{ {\bf e_{n+1}}, \ldots, {\bf e_N} \}$ of the orthogonal 
complement of  $\mbox{span}(\{ {\bf e_a} \})$ in $V$. Even though our attempts to 
extend ${\cal P}_q$ with functions
failed, now we are extending with $N$-tuples.
Define the obvious projection operator 
 $\pi : V \to {\cal P}_q$ 
 by
 $\pi {\bf e_1} =  e_1, \ldots, \pi {\bf e_n} =  e_n$, 
 $\pi {\bf e_{n+1}} = \ldots = \pi {\bf e_N} = 0 \in {\cal P}_q$;
clearly $\pi \iota = I$, the identity operator on ${\cal P}_q$.

Recall that the essential step in the proof of theorem 5 is extension of 
$Q$ to an $N \times N$ orthonormal
matrix $\tilde{Q}$ by appending any $N - n$ orthonormal rows. In this way 
$\tilde{A}_1,\ldots,\tilde{A}_d$, 
$N \times N$ symmetric commuting extensions of the $A_i$ are constructed 
in theorem 9, where $\tilde{A}_i = \tilde{Q} \Lambda_i \tilde{Q}^t$. 
Since the $\tilde{A}_1,\ldots,\tilde{A}_d$ mutually commute and are symmetric, there
exist $N$ orthonormal common eigenvectors 
${\bf u}_\alpha \in V$, such 
that $\tilde{A}_i {\bf u}_\alpha = \Lambda_{i \alpha } {\bf u}_\alpha$. 
The matrices $\tilde{A}_i$ are given in the basis $\{{\bf e_\alpha} \}$ and 
$\tilde{Q}$ is the transformation
between this basis and the eigenvector basis $\{ {\bf u}_\alpha \}$.
Note that the rows of $\tilde{Q}$ give the coordinates of the vectors 
$\{ {\bf e}_{\alpha} \}$
in the basis $\{ {\bf u}_{\alpha} \}$, hence the extension of 
$Q$ to $\tilde{Q}$ by adding
arbitrary $N-n$ orthonormal rows is nothing but the extension of 
$\{ {\bf e_a} \}$ to $\{ {\bf e_\alpha} \}$ 
the orthonormal basis of $V$ described above.
The reader is reminded that at present we are
assuming the existence of a cubature formula hence the 
eigenvalues $\Lambda_{i \alpha } $ are defined in (\ref{EigXconnect}).
We can now state:

\noindent{\bf Theorem 10}. The commuting extensions 
 of theorem 9, $\{ \tilde{A}_i \}$,  satisfy the following
compatibility condition with the operators $\{ \chi_i \}$,
\be  \tilde{A}_i \iota p = \iota x_i p = \iota \chi_i p  ,
    \qquad \forall p \in {\cal P}_{q-1}\ , \la{exstr}\ee
where ${\cal P}_{q-1}$ is 
regarded in the natural way as a subspace of ${\cal P}_q$.

\noindent{\bf Note:} Applying $\pi$ to \r{exstr} gives
$\pi \tilde{A}_i \iota p = \chi_i p $, for $p \in{\cal P}_{q-1}$, which 
is automatic as $\tilde{A}_i$ is an extension of $A_i$. However, \r{exstr} 
contains more information than this, and is not true for an arbitrary extension
of $A_i$. 

\noindent{\bf Proof.} We first prove that for any $p \in {\cal P}_q$ and any eigenvector ${\bf u}_{\alpha}$,
 \be \langle \iota p | {\bf u}_\alpha \rangle =  \sqrt{w_\alpha} p(x_\alpha) \ , \label{ProjLemm2a} \ee
where $w_\alpha$, $x_\alpha$, are the weights and nodes of the cubature formula whose existence is assumed 
in theorem 9. 
We already noted that the rows of $\tilde{Q}$ from 
the proof of theorem 9 give the coordinates
of the $N$ vectors ${\bf e}_\alpha$ in the basis $\{ {\bf u}_\alpha \}$, 
in particular the rows of $Q$ give the coordinates of ${\bf e_a} = \iota e_a$, $a = 1,\ldots,n$. 
Recall that $Q_{a\alpha} = \sqrt{w_\alpha} e_a(x_\alpha)$,
thus \r{ProjLemm2a} is proven for basis elements $e_a \in {\cal P}_q$. The proof for general 
$p \in {\cal P}_q$ follows by linearity. 

We now expand the left hand side of (\ref{exstr}) in the basis $\{ {\bf u_\alpha} \}$.
For any $p \in {\cal P}_{q-1}$
\be \langle \tilde{A}_i \iota p | {\bf u_\alpha}   \rangle    =
   \Lambda_{i \alpha} \langle  \iota p | {\bf u_\alpha}   \rangle =
   \Lambda_{i \alpha} \sqrt{w_\alpha} p(x_\alpha) =  \sqrt{w_\alpha} (x_\alpha)_i p(x_\alpha)  \ , \ee
using (\ref{ProjLemm2a}) and (\ref{EigXconnect}) respectively in the last two steps. 
To expand the right hand side of (\ref{exstr}) note that $ x_i p = \chi_i p \in {\cal P}_q  $ for all
$p \in {\cal P}_{q-1}$. Invoking (\ref{ProjLemm2a}) again we obtain
\be
   \langle \iota x_i p | {\bf u_\alpha} \rangle = \sqrt{w_\alpha} (x_\alpha)_i p(x_\alpha) \ ,
\ee
which completes the proof.
$\bullet$

\noindent
Note that taking $p = 1$ in (\ref{ProjLemm2a}) gives
$w_\alpha = \langle \iota 1 | {\bf u}_{\alpha} \rangle^2$. In theorem 9 we saw that the eigenvalues
of the commuting extensions are related to cubature nodes, here we obtain a relation between
the weights and common eigenvectors.

We shall see in section 5.3 that with an appropriate choice of basis the compatibility condition 
of theorem 10
implies that the commuting extensions have certain off-diagonal zero blocks;
in particular this special structure
aids computation of commuting extensions.

The obvious question to ask at this stage is whether there is a converse
of theorems 9 and 10.  That is, suppose we have 
$\tilde{A}_1,\ldots, \tilde{A}_d$, $N \times N$ symmetric commuting extensions
of $A_1,\ldots, A_d$, which satisfy the compatibility condition
 $  \tilde{A}_i \iota p = \iota \chi_i p = \iota x_i p ,
    \ \forall p \in {\cal P}_{q-1}\ $.
Can we build a cubature rule, without {\em a priori} assumption of its existence,
using the eigenvalues and eigenvectors of $\tilde{A}_1,\ldots, \tilde{A}_d$?
In theorem 11 we give an affirmative answer to this. The treatment follows the 
presentation on Gaussian
quadrature from section 5.1; in particular we start with a $\delta$ lemma. 
Note that given the commuting matrices $\tilde{A}_i$ we can find their diagonal representations $\Lambda_i$, 
but we do not assume
in advance any connection of the $\Lambda_i$ 
with cubature nodes.

\noindent{\bf Lemma 3} (The multidimensional $\delta$ lemma)
Suppose the commuting extensions satisfy the compatibility condition 
$\tilde{A}_i \iota \grave{p} = \iota x_i \grave{p} = \iota \chi_i \grave{p}$
for all $i = 1, \ldots, d,$ and for all $\grave{p} \in {\cal P}_{q-1}$. Then,
for any $p \in {\cal P}_q$ 
\be \langle \iota p |   {\bf u}_\alpha \rangle  =  \langle \iota 1 |  {\bf u}_\alpha \rangle  \ p( \lambda_\alpha) \ ,  \label{ExtProjForm}
\ee
where the points
$\lambda_\alpha\in{\bf R}^d$ have entries $(\Lambda_{1 \alpha },\ldots,\Lambda_{d \alpha })$, all eigenvalues
of $\tilde{A}_1, \ldots, \tilde{A}_d $,
satisfying
$A_i {\bf u_\alpha} = \Lambda_{i \alpha } {\bf u_\alpha}$. 

\noindent{\bf Proof.}
We prove (\ref{ExtProjForm}) for monomials 
$p = x_1^{m_1}x_2^{m_2}\ldots x_d^{m_d} $. 
For $p = 1$ the statement is trivial. For any other monomial 
$p \in {\cal P}_q$ we can
write $p = x_i \grave{p}$, where 
$\grave{p} \in {\cal P}_{q-1}$. Then
\be
 \langle {\iota p} | {\bf u}_\alpha \rangle = \langle \iota x_i \grave{p} | {\bf u}_\alpha \rangle = 
 \langle \tilde{A}_i \iota \grave{p} | {\bf u}_\alpha \rangle =  
 \langle  \iota \grave{p} | \tilde{A}_i {\bf u}_\alpha \rangle
 = \Lambda_{i \alpha } \  \langle {\iota \grave{p}} |   {\bf u}_\alpha \rangle \ .  \label{RecStep}
\ee
Here the compatibility condition was used in the second step. Repeated application of (\ref{RecStep}) completes 
the proof for monomial $p$; the full result follows by linearity. $\bullet$

\noindent{\bf Note:}
Recall that we do not know how to relate $V$ to a space of 
polynomials (or other functions) in a way
which gives commuting extensions of the operators $\chi_1,\ldots,\chi_d$.
In particular,
we can not interpret the eigenvectors ${\bf u_\alpha}$ as polynomials (or other functions).
Thus our present state of understanding allows us to view the ${\bf u_\alpha}$ as
``$\delta$ {\em vectors}'' in $V$ and not as $\delta$ {\em functions}, which was possible
in the 1-dimensional case. Moreover we can not identify ${\cal P}_{q+1}$ with a subspace of $V$.
Hence,
in contrast to the one dimensional case, we restrict $p \in {\cal P}_q$ in the multidimensional $\delta$ 
lemma thereby losing the immediate 
connection between cubature nodes and roots of polynomials in ${\cal P}_q^\perp $.

We are now fully prepared for the converse statement to theorems 9 and 10:

\noindent{\bf Theorem 11}. 
Let $A_1,\ldots,A_d$ be the representation of the 
operators $\chi_1,\ldots,\chi_d$ in an orthonormal
basis $\{e_a\}$ of ${\cal P}_q$. 
Let $\tilde{A}_1,\ldots,\tilde{A}_d$ be 
$N\times N$ symmetric commuting extensions of
$A_1,\ldots,A_d$ satisfying the compatibility condition
$ \tilde{A}_i \iota p = \iota x_i p = \iota \chi_i p$  
for all $p  \in {\cal P}_{q-1}$. 
Then for every polynomial $f$ in ${\cal P}_{2q+1}$ 
\be \int_\Omega w(x) f(x) d^dx =
  \sum_{\alpha=1}^N \langle \iota 1 |{\bf u_\alpha}\rangle^2  f(\lambda_\alpha) \ . 
   \la{QuadRule}\ee 
Here the ${\bf u_\alpha}$ are joint eigenvectors of $\tilde{A}_1, \ldots, \tilde{A}_d $,
satisfying
$A_i {\bf u_\alpha} = \Lambda_{i \alpha } {\bf u_\alpha}$, 
and 
the points
$\lambda_\alpha\in{\bf R}^d$ have entries $(\Lambda_{1 \alpha },\ldots,\Lambda_{d \alpha })$.

\noindent{\bf Proof}. Recall the symmetric bilinear forms $X_i$ on ${\cal P}_q$ 
defined in \r{qfd}. Given the commuting extensions, we can
introduce symmetric bilinear forms $\tilde{X}_i$ on 
$V$ defined by $\tilde{X}_i({\bf u}, {\bf v})= \langle {\bf u} | \tilde{A}_i {\bf v} \rangle$. 
Since the $\tilde{A}_i$ 
are extensions of the $A_i$ we have  
$X(p_1,p_2)=\tilde{X}_i(\iota p_1, \iota p_2 )$ for all $p_1,p_2$ in ${\cal P}_q$.
The $\tilde{X}_i$ are simultaneously diagonalized in the basis $\{ {\bf u_\alpha } \}$, 
$\tilde{X}_i({\bf u_\alpha},{\bf u_\beta})=\Lambda_{i \alpha }\delta_{\alpha\beta}$.

It is sufficient to prove the statement of the theorem for monomials. 
For $f = 1$ we have 
\be
\int_\Omega w(x)   d^dx  = 
\langle 1|1\rangle = \langle \iota 1 | \iota 1 \rangle 
= \left\langle  \sum_{\alpha=1}^N \langle \iota 1  | {\bf u_\alpha} \rangle {\bf u_\alpha} \left|
         \sum_{\beta=1}^N \langle \iota 1  | {\bf u_\beta} 
     \rangle {\bf u_\beta} \right.\right\rangle
= \sum_{\alpha=1}^N  \langle \iota 1 | {\bf u_\alpha} \rangle^2  \ ,\la{77}\ee
by the orthonormality of the ${\bf u_\alpha}$. Note that in the second
expression in \r{77} the inner product is taken in ${\cal P}_q$, in
subsequent expressions it is taken in $V$.

Any other monomial in ${\cal P}_{2q+1}$ can be written
in the form $f = x_i f_1  f_2 $ for some monomials $f_1,f_2\in {\cal P}_q$ and some
$i$. Note that use of the multidimensional $\delta$ lemma is possible since we assume the
$\tilde{A}_i$ satisfy the compatibility condition, and so 
\bea
\int_\Omega w(x) f(x)  d^dx  
&=& \int_\Omega w(x) f_1(x) x_i f_2(x)\ d^dx     
 = X_i ( f_1, f_2 ) \nonumber\\
&=& \tilde{X}_i ( \iota f_1 , \iota f_2 )  
 =  \tilde{X}_i \left(  \sum_{\alpha=1}^N \langle \iota f_1 | {\bf u_\alpha}\rangle {\bf u_\alpha} ,
             \sum_{\beta=1}^N \langle \iota f_2 | {\bf u_\beta}\rangle {\bf u_\beta} 
     \right)   \nonumber \\
&=& \sum_{\alpha=1}^N  \Lambda_{i \alpha }
                        \langle \iota f_1 |{\bf u_\alpha}\rangle 
                       \langle \iota f_2 | {\bf u_\alpha}\rangle 
 = \sum_{\alpha=1}^N  \langle \iota 1 | {\bf u_\alpha}\rangle^2  \Lambda_{i \alpha }
                         f_1(\lambda_{\alpha }) 
                        f_2(\lambda_{\alpha }) \nonumber\\
&=& \sum_{\alpha=1}^N  \langle \iota 1 | {\bf u_\alpha} \rangle^2 
                         f(\lambda_{\alpha }) 
                        \ .   
\eea
Thus \r{QuadRule} is proven for monomial $f$; the full result follows by linearity.  
$\bullet$

Theorems 9,10,11 give the main result of this paper, that $N$ point, odd order
cubature formulae with positive weights are equivalent to symmetric commuting extensions, 
satisfying the compatibility condition, of matrix representations of the operators
$\chi_1,\ldots,\chi_d$.

\subsection{Discussion and Consequences}

Our findings give a new computational approach to the derivation of cubature 
formulae. If appropriate commuting extensions
are numerically found their simultaneous diagonalisation will give the cubature rule
in (\ref{QuadRule}).
To numerically obtain the matrices $\{ A_i \}$ we introduce an orthonormal 
basis of ${\cal P}_q$ consisting of an orthonormal basis of ${\cal P}_0$ 
(a constant function $e_1$ with $\| e_1 \| = \sqrt{\langle e_1|e_1 \rangle} = 1$), 
extended to one of ${\cal P}_1$,
extended to one of ${\cal P}_2$ etc., $i.e.$ a  basis $\{ e_a \}$, 
$a=1,\ldots,n=\pmatrix{d+q \cr d\cr}$,  
such that 
\be
\begin{array}{l}
{e_1}~ \mbox{is an orthonormal basis of}~{\cal P}_0  \\
{e_1,\ldots,e_{d+1}}~ \mbox{is an orthonormal basis of}~{\cal P}_1  \\
{e_1,\ldots,e_{\frac12(d+1)(d+2)}}~ \mbox{is an orthonormal basis of}~{\cal P}_2  \\
{\rm etc.}
\end{array}
\ee
Such a basis can be obtained from the monomials 
$\{ x_1^{m_1} \ldots x_d^{m_d}  \}, \ m_1 + \ldots + m_d \leq q $, by the Gram-Schmidt 
procedure.
Note that all basis elements $e_a$ of degree $m$ or more are orthogonal to 
${\cal P}_{m-1}$.
This choice of basis and the fact that $\{ A_i \}$ represent the 
operators $\{ \chi_i \}$
imply that the $ A_i $ have tridiagonal block form as in \r{trid}, 
with $q+1$ blocks on the 
diagonal.
In the notation of section 4, $n_1 = 1$ and for $m = 2, 3, \ldots, q+1$, 
$n_m = \mbox{dim}{\cal P}_{m-1} - \mbox{dim}{\cal P}_{m-2}$
($n_1=1$, $n_2=d$, $n_3 = d(d+1)/2$ etc). 
Moreover, the commutator of any pair of 
$A_i$'s is zero apart from a single block of size 
$\pmatrix{q+d-1 \cr d-1\cr}\times\pmatrix{q+d-1 \cr d-1\cr}$ 
in the bottom right hand corner
(note $\pmatrix{q+d-1 \cr d-1\cr}= n_{q+1} = \mbox{dim}{\cal P}_q - \mbox{dim}{\cal P}_{q-1}$).  
The compatibility condition, together with 
the fact that ${\bf e_{n+1}}, \ldots, {\bf e_{N}}$ are orthogonal to $\iota {\cal P}_q$, imply that the bottom 
left $(N-n) \times {\rm dim}({\cal P}_{q-1})$ block of 
$\tilde{A}_i$ is zero and by symmetry so is the corresponding block in the upper 
right. Thus the commuting
extensions $\{ \tilde{A}_i \}$ we seek are
precisely those in tridiagonal block form as in section 4. 
Note also that since our first basis element $e_1$ is a constant polynomial the 
cubature weights are obtained from the first
entries of the eigenvectors ${\bf u_\alpha}$, 
$w_\alpha = \langle \iota 1 | {\bf u_\alpha} \rangle^2 =  \frac{1}{e_1^2} 
\langle \iota e_1 | {\bf u_\alpha} \rangle^2 = 
\left( \frac{ ({\bf u_\alpha})_1 }{e_1} \right)^2 $.

In the case $d=2$, we note that the matrices $A_1,A_2$ have 
$n_{r-1}=q$ and  $n_r=q+1$ in the notation of section 4, and thus 
from \r{nbd}, which is based on counting degrees of freedom, 
the expected size of the commuting extensions is
\be
N \ge n +\frac{q(q+1)}{6} \ . 
\ee
Using $n={\rm dim}{\cal P}_q=\frac12(q+1)(q+2)$ we obtain
\be
N \ge \frac{(2q+2)(2q+3)}{6} = \frac13 {\rm dim}{\cal P}_{2q+1} \ . 
\la{75}\ee
This is exactly the number of nodes we expect from counting degrees of freedom in a
2-dimensional cubature formula of degree $2q+1$ . 

For $d > 2$ the $A_i$ have a more subtle 
structure that requires refinement of the discussion leading to equation \r{nbd}.
However, the equivalence of cubature formulae and commuting extensions (satisfying the
compatibility condition) allows us to estimate $N$ by easily counting degrees of freedom
in a general d-dimensional cubature formula of degree $2 q + 1$. Thus,
\be N \geq \frac{1}{d+1}\mbox{dim}{\cal P}_{2q+1}     \ .\ee 
We emphasize again that such calculations are not rigorous and the inequalities 
obtained in this
way can serve only as recomendations for choice of $N$, indeed certain cubature 
formulas with less points are known. 

In section 6 we shall give some first examples of computation of 
cubature nodes using our approach. But before we do this we 
present two theoretical consequences of
the equivalence between cubature formulae and 
commuting extensions. 

\noindent{\bf Theorem 12}. Let $N$ be the number of nodes in a degree $2q+1$, 
$d$-dimensional positive weight cubature rule. Then
\be N \geq \pmatrix{d+q \cr d\cr}
 + \frac12 {\rm max}_{i,j} {\rm rank}([A_i,A_j])\  ,
   \la{lub}\ee 
where $A_1,\ldots,A_d$ are the matrix representations of 
the operators $\chi_1,\ldots,\chi_d$ on ${\cal P}_q$.  

\noindent{\bf Proof}. By theorem 9 an $N$ point  cubature rule gives  $N\times N$ 
commuting extensions of the matrices $A_i$. By theorem 2, section 2, the size
of such 
extensions is at least ${\rm dim}{\cal P}_q+\frac12 
{\rm max}_{i,j} {\rm rank}([A_i,A_j])$. $\bullet$  

\noindent{\bf Notes:} (1) As mentioned in the introduction, theorem 12 has
its origins in the work of M\"oller \c{Moller}. A statement of the 
result in a form that clearly corresponds to our statement can be found 
in \c{ex1}, which cites \c{ex2} and \c{Xu}.
Our proof, however, is a substantial
simplification. 
(2) It is informative to compare the lower bound of theorem 12
with estimates based on parameter counting.
As a consequence of our previous remarks on the block structure of $[A_i, A_j]$
\be{\rm rank}([A_i,A_j])\le \pmatrix{d-1+q \cr d-1\cr} =
   \frac{d}{d+q} \pmatrix{d+q \cr d\cr}\ ,       \ee
so the second term in \r{lub} is typically a small fraction of the first term.
Consequently for large $q$ the right hand side of \r{lub} is much  smaller than 
the number of nodes we expect from counting degrees of freedom in a degree $2q+1$, 
$d$-dimensional cubature formula, which is
\be 
\left\lceil \frac1{d+1}\pmatrix{d+2q+1 \cr d\cr}\ \right\rceil.
\ee
This comparison indicates why the lower bound on the number of points needed for a
cubature formula is rarely attained.

\noindent{\bf Theorem 13}. 
Let $A_1,\ldots,A_d$ be matrix representations of $\chi_1,\ldots,\chi_d$.  
In any degree $2q+1$, $d$-dimensional, positive weight cubature rule, 
and for each $i$, there is a node $x_\alpha$ with $(x_\alpha)_i$
less than or equal to the smallest eigenvalue of $A_i$,  
and a node $x_\beta$ with 
$(x_\beta)_i$ greater than or equal to the largest eigenvalue
of $A_i$. 

\noindent{\bf Proof}. By theorem 9 a cubature rule of degree $2q+1$ 
gives commuting extensions of the matrices $A_i$ with the nodes composed of
the eigenvalues of the extended matrices. By theorem 6 in section 2 the 
smallest/largest eigenvalue of the extended matrices is less/greater than
or equal to the smallest/largest eigenvalue of the matrices before extension. 
$\bullet$

\noindent{\bf Note:} As far as we are aware this theorem is not even known 
for $d=1$. For $d=1$ the theorem says that any $N$-point, positive weight,
degree $2q+1$ quadrature rule must have a node less/greater than or equal to the 
smallest/largest Gaussian quadrature node. Thus Gaussian quadrature has the 
property that the span of the nodes is the smallest possible, amongst all 
positive weight quadrature rules, with any number of points, that are exact 
to the same degree.

\section{Examples}

In this section we briefly discuss 2 examples of computing cubature rules
via commuting extensions, both are in 2 dimensions. In section 6.1 we consider 
the classic question, first studied by Radon \c{radon}, of finding 7 point 
cubature rules which are exact for polynomials up to degree 5 (${\rm dim}
{\cal P}_5=21$
for $d=2$). This involves $7\times 7$ extensions of a pair of $6\times 6$
matrices of the tridiagonal block form of section 4, and we provide a reliable algorithm for this. 
In section 6.2 we present some new 
cubature rules for integration on the entire
plane with weight function $w(x,y)=e^{-x^2-y^2}$. These were computed 
by commuting extension techniques, though, as explained in section 3, our
current algorithms for this are poor, so we give limited details.  

\subsection{Radon type formulae}

Given a region $\Omega$ on the plane and a suitable weight function
$w(x_1,x_2)$, a Radon formula is a cubature rule of the form 
\be
\int_\Omega  w(x_1,x_2) f(x_1,x_2)\  d^2x
= \sum_{\alpha=1}^7 w_\alpha f(x_{1\alpha},x_{2\alpha}) 
\ee
which is exact for all polynomials of degree no more than 5. 

To construct such formulae we proceed as follows: First choose a 
basis $e_1,e_2,e_3,e_4,e_5,e_6$ of ${\cal P}_2$ which is orthonormal 
with respect to the inner product 
\be
\langle a | b \rangle =
\int_\Omega w(x_1,x_2)a(x_1,x_2)b(x_1,x_2)\  d^2x\ ,
\ee
and such that $e_1$ is of degree $0$, $e_2,e_3$ are of degree $1$, and 
$e_4,e_5,e_6$ are of degree $2$. Typically we do this by applying 
Gram--Schmidt orthonormalization to the basis $1,x_1,x_2,x_1^2,x_1x_2,x_2^2$.

Having constructed the orthonormal basis we compute the matrices 
$A_1$ and $A_2$
via:
\bea 
(A_1)_{ij} &=& \int_\Omega w(x_1,x_2) e_i(x_1,x_2) x_1 e_j(x_1,x_2) \  d^2x\ ,
\la{defA1}\\
(A_2)_{ij} &=& \int_\Omega w(x_1,x_2) e_i(x_1,x_2) x_2 e_j(x_1,x_2) \  d^2x\ .
\la{defA2}\eea
By orthogonality we will have $(A_1)_{1i}=(A_1)_{i1}=(A_2)_{1i}=(A_2)_{i1}=0$
for $i=4,5,6$ (this is the meaning of ``tridiagonal block form'' in this case),
and the commutator $C=[A_1,A_2]$ will be all zero except for a single 
antisymmetric $3\times 3$ block in the lower right hand corner. 
The commuting extensions $\tilde{A_1},\tilde{A_2}$ should  take the form 
given in \r{fsce}, 
where $a,b$ are 6-dimensional column vectors, with the first 3 entries vanishing
(these are the ${\rm dim}{\cal P}_{q-1} \times (N-n)$ zero blocks described in
section 5.3), 
and $\alpha,\beta$ scalars.
Following the arguments in lemma 1 from section 2,
$a,b,\alpha,\beta$ must satisfy \r{re1}-\r{re2}.
To solve these equations we 
proceed exactly as in the lemma.
First we construct a specific pair of 6-dimensional vectors $v,w$ satisfying 
\be [A_1,A_2]+vw^T-wv^T = 0\ . \ee
Because of the special structure of $[A_1,A_2]$, $v,w$ can be taken to have zeros
in their first 3 entries, and determining the other entries is equivalent to
finding 2 vectors in ${\bf R}^3$ with a given cross product. 
Clearly $v,w$ are linearly independent and give a basis for ${\rm Im}([A_1,A_2])$.
Once such  $v,w$ have been found, \r{re1} gives
\be a=\lambda v + \mu w \ , \qquad 
    b=\nu v + \rho w \ , \la{defab}\ee
where $\lambda \rho - \mu \nu = 1 $, and  \r{re2} gives the requirement 
\be
 (\beta \lambda - \alpha \nu) v 
+ (\beta \mu  - \alpha  \rho) w 
+ \nu A_1 v 
+ \rho A_1 w 
- \lambda A_2 v  
- \mu A_2 w 
= 0 \ ,
\ee
or
\be
\pmatrix{ v & w & A_1 v & A_1 w & A_2 v & A_2 w \cr}
\pmatrix{ \beta \lambda - \alpha \nu \cr
          \beta \mu  - \alpha  \rho  \cr
         \nu \cr
         \rho  \cr
      - \lambda  \cr
      - \mu \cr }
= 0 \ .
\la{fr1}\ee
Since $v,w$ have zeros in their first 3 entries, and $(A_1)_{1i}=(A_2)_{1i}=0$
for $i=4,5,6$, the entire first row of the matrix 
$\pmatrix{ v & w & A_1 v & A_1 w & A_2 v & A_2 w \cr}$ is zero. So
it is singular and a nontrivial solution of 
$\pmatrix{ v & w & A_1 v & A_1 w & A_2 v & A_2 w \cr} k= 0 $ is 
guaranteed. To complete construction of
a commuting extension all we need is to find
$\alpha,\beta,\lambda,\mu,\nu,\rho$ with $\lambda \rho - \mu \nu = 1 $,
and such that the column vector in \r{fr1}
is in the kernel of the matrix in \r{fr1}. Using the freedom to
rescale vectors in the kernel, it is straightforward to conclude
that given a nontrivial vector $k$ in the kernel of the matrix,
there will be an associated commuting extension if and only if
$k_3k_6-k_4k_5>0$, and we must take
\bea 
&\lambda=-ck_5\ , \qquad 
\mu=-ck_6\ , \qquad 
\nu=ck_3\ , \qquad 
\rho=ck_4\ ,  &\nonumber\\
&\alpha=c^2(k_2k_5-k_1k_6)\ , \qquad 
\beta=c^2(k_1k_4-k_2k_3)\ , &\la{defall}\\
&{\rm where} \quad c=\frac1{\sqrt{k_3k_6-k_4k_5}}&
\nonumber
\eea

To summarize: to find a Radon formula one should:
(1) Construct the orthonormal basis $\{e_a\}$. 
(2) Construct the matrices $A_1,A_2$.
(3) Find $6$-dimensional vectors $v,w$ with top 3 entries zero,
and such that $[A_1,A_2]+vw^T-wv^T = 0$.
(4) Compute the kernel of the matrix $\pmatrix{ v & w & A_1 v & A_1 w & A_2 v & A_2 w \cr}$. 
Then for each vector $k$ in the kernel with $k_3k_6-k_4k_5>0$ there is an 
associated commuting extension given by \r{fsce}, where 
$a,b$ are given by \r{defab} and 
$\alpha,\beta,\lambda,\mu,\nu,\rho$ by \r{defall}. The commuting extensions 
can be simultaneously diagonalized, using the algorithm in \c{c}, 
to obtain the nodes and weights of the cubature rule.

This procedure should be compared with Radon's original work \c{radon}. 
In practice we have not found  a case in which there is a vector in the 
kernel with $k_3k_6-k_4k_5\le 0$, but have no explanation why this is so. 
Generically the kernel is one dimensional 
and there is a unique Radon formula.  For the case of 
$\Omega$ equal to the circle or the 
square each with uniform weight $w(x_1,x_2)=1$,  
the kernel has dimension 2, giving a 
one parameter family of Radon formulas. In the case of the square
with vertices $(-1,-1),(-1,1),(1,1),(1,-1)$, 
all the Radon formulae 
have a single node at the origin with weight $\frac87$ and 3 
pairs of diametrically opposed nodes on the circle $x_1^2+x_2^2=\frac{14}{15}$. 
(In the literature there is often only mention of the case 
of Radon formulae for the square when one
pair of nodes lies on one of the coordinate axes.)
In the case of the unit circle, there is
a single node at the origin with weight $\frac\pi4$ and 3 
pairs of diametrically opposed nodes on the circle $x_1^2+x_2^2=\frac23$\ ;
in this case there is full rotational symmetry, and the different formulae
are related by rotation. 


As an explicit example of nonstandard 
Radon formulae, we give here some results for 
the square with vertices $(-1,-1),(-1,1),(1,1),(1,-1)$,
with a square with vertices $(\frac25-r,\frac35-r),(\frac25-r,\frac35+r),
(\frac25+r,\frac35+r),(\frac25+r,\frac35-r)$ removed. 
Here $0\le r\le \frac25$ and the weight is uniform. Figure 1 displays
the location of the nodes for $r=\frac{i}{20}$, $i=0,1,\ldots,8$, along 
with the boundaries of the squares removed, and the circle 
$x_1^2+x_2^2=\frac{14}{15}$ on which the nodes lie for $r=0$. For $r=0$ we
in fact plot the nodes for a very small nonzero value of $r$, so as to have a 
unique result. For values of $r$ up to $\frac14$ all the nodes are inside
the domain of integration. For $r=\frac3{10}$ and $r=\frac7{20}$ 
one node lies inside the small square that is removed to obtain the 
cubature domain. For $r=\frac25$ one node lies outside the larger square,
at $(x_1,x_2)\approx(0.1844,1.0360)$, but has low weight (about $3.25\%$ of the 
sum of the weights). 

\begin{figure}
\centerline{\includegraphics[width=15cm]{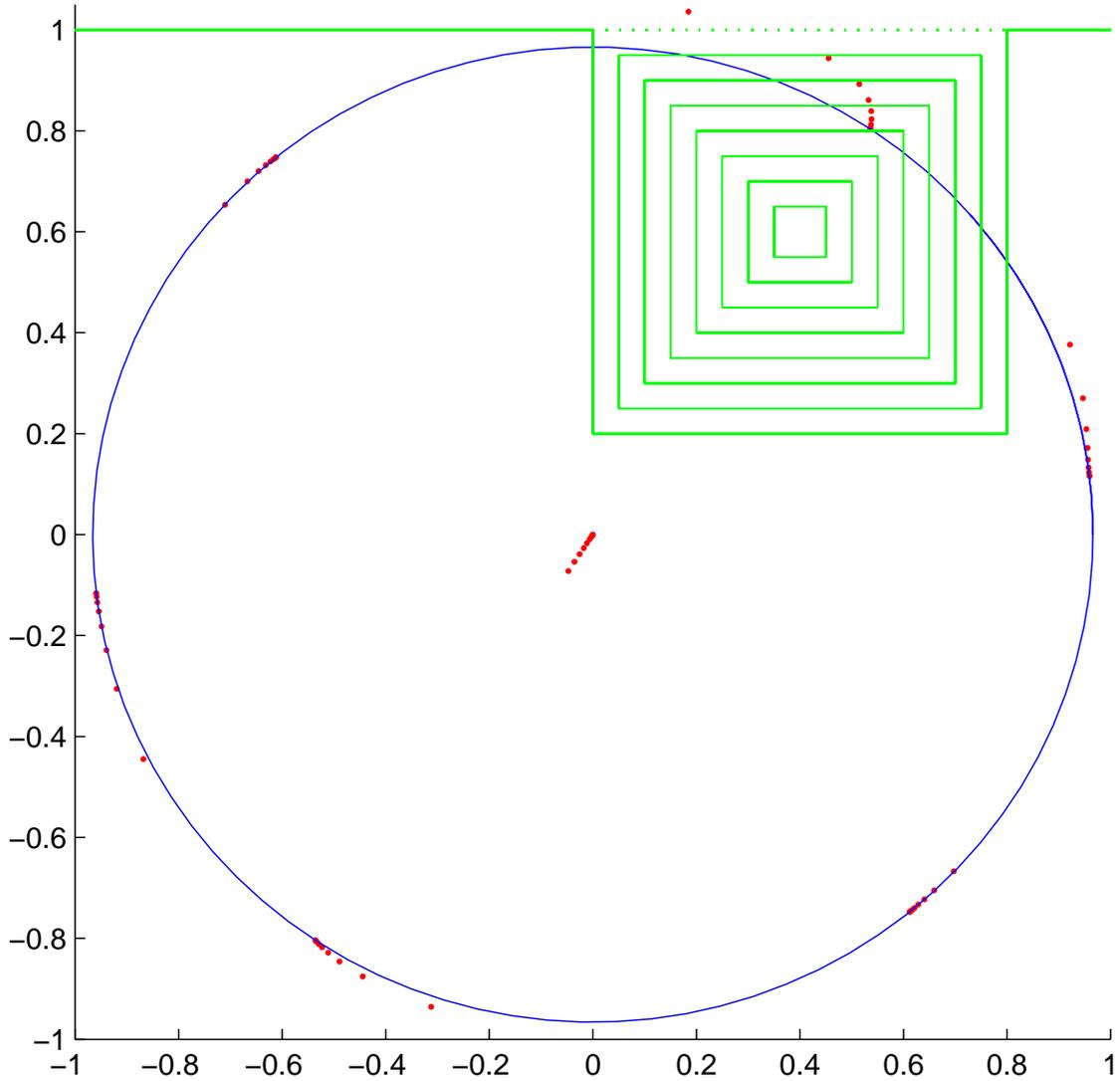}}
\caption{Radon formulae on a square with a smaller square removed, with uniform weight.
As the size of the smaller square tends to zero, one node tends to the origin and 6 to 
the circle $x_1^2+x_2^2=\frac{14}{15}$. As the size of the smaller square is increased, one
node moves outside the domain of integration, and eventually out of the larger square.}
\end{figure}

Recall that our approach for constructing cubature formulae 
involves first constructing commuting
extensions, and then diagonalizing them. For the rational values of $r$ discussed in
the previous paragraph, the commuting extensions can be constructed exactly, and 
we do not need to apply the numerical algorithms for commuting extensions 
from section 3 (the $S(Q,\Lambda_i)$ algorithm does not work particularly well for the matrices in
this subsection). Even though the commuting extensions can be found exactly, 
finding the cubature nodes (the eigenvalues of the extended matrices) can only be
done numerically.

\subsection{Gaussian weight on the plane}

As already explained, 
to construct degree $2q+1$ cubature formulae on a region $\Omega$ in the 
plane, we need to find symmetric commuting extensions (satisfying the necessary
compatability conditions) of the  $n\times n$  matrices defined in 
\r{defA1}-\r{defA2}, where 
$n={\rm dim}{\cal P}_q=(q+1)\left(\frac12q+1\right)$. By counting degrees of freedom,
we expect $N$-dimensional symmetric extensions if 
$N\ge (q+1)\left(\frac23q+1\right)$, see the discussion leading to \r{75}.
We work in a basis of the 
type discussed in subsection 5.3, so the matrices are all of the tridiagonal
block form of section 4. Using the notation of section 4, 
we extend each matrix by adding 2 blocks, $a_i$  of 
size $(q+1)\times(N-n)$, and $\alpha_i$ (symmetric), 
of size $(N-n)\times(N-n)$, here $i=1,2$. 
Thanks to the freedom in choice of extensions described in theorem 3, one
of the symmetric added blocks $\alpha_i$ can be chosen diagonal. Taking this
into account we expect 
$N-n\ge\frac16q(q+1)$. Even so, the number of entries added 
in the blocks $a_1,a_2,\alpha_1,\alpha_2$ grows
as $q^4$. 

For applications in quantum mechanics, cubature on the plane 
$\Omega ={\bf R}^2$ with 
weight function $w(x_1,x_2)=e^{-x_1^2-x_2^2}$ is important. For this case
we have computed degree
11, 13, 15 and 17 cubature formulas with 26, 35, 46 and 57 nodes respectively,
by simultaneous diagonalization of commuting extensions. These were found using
the gradient flow approach mentioned at the begining of section 3.
For degrees 11 and 13 we have also 
succeeded to compute the necessary commuting extensions using the $S(Q,\Lambda_i)$ algorithm on which 
section 3 focused. To take into account the zero blocks of 
the commuting extensions an extra  term 
\be \sum_{i=1}^d \sum_{a=1}^{{\rm dim}{\cal P}_{q-1}}
    \sum_{b=n+1}^N \left( (\tilde{Q} \Lambda_i \tilde{Q}^T)_{ab}
    \right) ^2  \ee
was added to $S$ of \r{Sdef}, this gave a substantial improvement in 
performance (compared to the same algorithm with $S$ not including the extra term). The nodes
and weights are available at \newline
\centerline{\tt http://www.math.biu.ac.il/$\sim$schiff/commext.html.} 
\newline Figure 2 displays the location of the 
nodes in the degree 15 and 17 formulae. The formulae found are
not necessarilly minimal, and are certainly not unique, because of rotational
symmetry. It seems our degree 13, 15 and 17 formulae are new, the one with 
order 17 having a smaller number of nodes compared to existing formulae of
the same order, see \c{enc}. 

\begin{figure}
\centerline{\includegraphics[width=8cm]{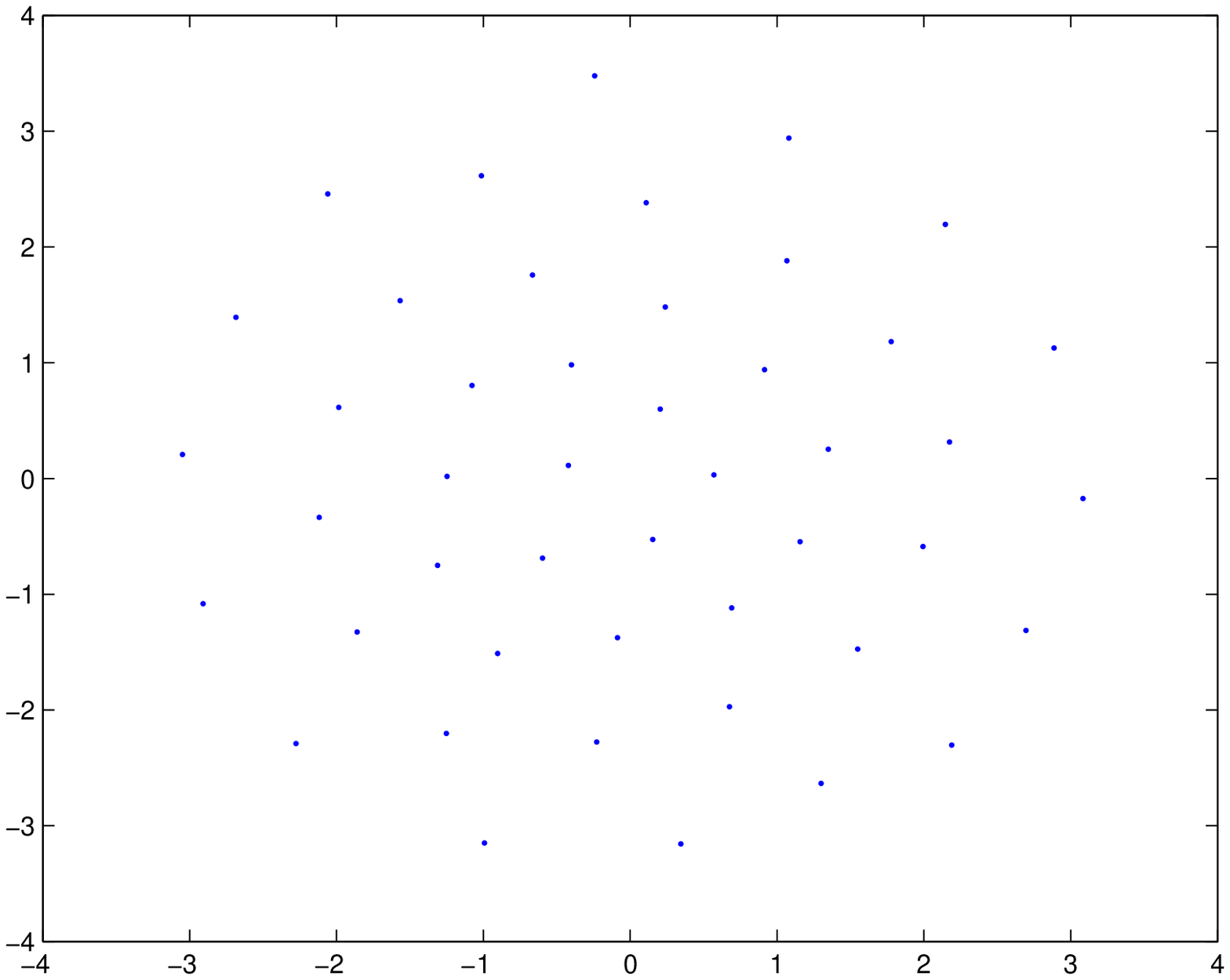}\includegraphics[width=8cm]{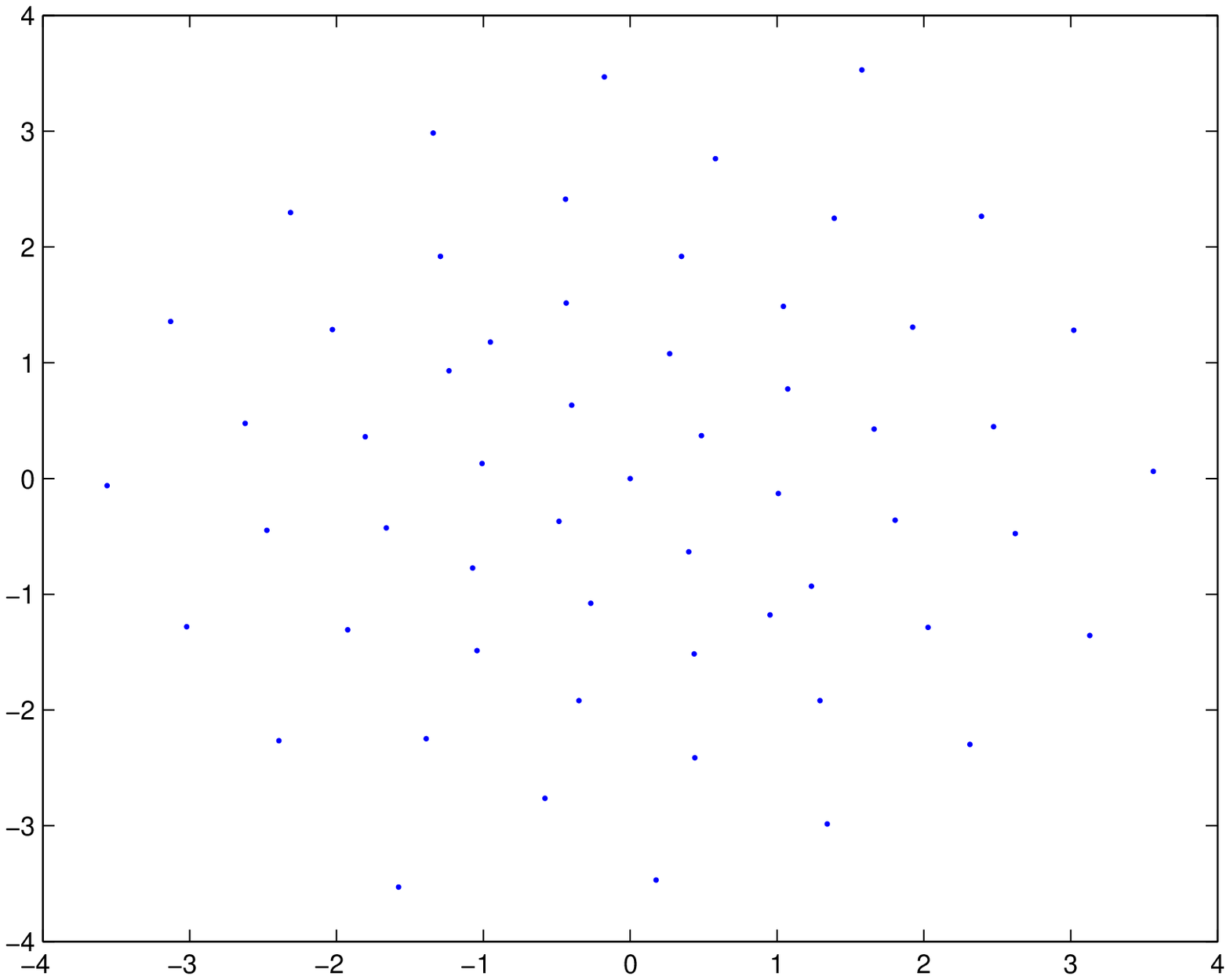}}
\caption{Nodes for a 46 point, degree 15 (left) and a 57 point, degree 17 (right)
cubature formula for the plane $\Omega={\bf R}^2$
with Gaussian weight $w(x_1,x_2)=e^{-x_1^2-x_2^2}$} 
\end{figure}

\section{Summary and open questions}

The central results of this paper are theorems 9,10,11, which prove the
equivalence of cubature formulae and commuting extensions satisfying the
compatibility condition (equivalent in an appropriate basis to
requiring certain zero blocks in the extension matrices). This raises the 
questions of existence and methods of computation for 
commuting extensions. Our knowledge of the
theory of commuting extensions is summarized by theorems 1 to 6, and in
section 3 we have described our initial attempts at their computation. 

There is clearly enormous potential for further work here. In the context
of our main topic, the connection between cubature formulae and commuting
extensions, there is one nagging question that we have indicated several
times in section 5: The vector space $V$ on which the commuting extensions
act does not yet have an interpretation as a space of functions (or maybe
even polynomials). For numerical work in quantum mechanics it would be a
major advantage if we could contstruct finite dimensional function spaces
containing the space of all polynomials of a given degree as a subspace,
on which the natural projections of the operators $x_i$ commute. The
existence (or nonexistence) of such spaces is a topic we hope to
investigate, see also \cite{cl}.

Another question left open in our work is that we have not given an 
existence proof of cubature formulae from the commuting extension viewpoint. 
Although theorem 1 in section 2 
guarantees the existence of commuting extensions of an arbitrary set of 
matrices, it does not guarantee extensions in the form we need to 
be able to apply theorem 11. An existence proof for commuting 
extensions of the required form would provide an alternative approach to
{\em Tchakaloff's theorem} \c{tchak}
that guarantees (for any suitable domain $\Omega$ and weight function $w(x)$) 
the existence of positive weight cubature rules that are exact for certain 
sets of functions. 

The numerical question of computing cubature formulae is now subsumed
under the more general question of computing commuting extensions,
likewise the open sore that there is no good way to predict the minimal number 
of points needed for a cubature rule is subsumed under the question of
finding the minimal dimension for commuting extensions. There are a 
number of points in the theory of commuting extensions which we feel may be
improved, for example, existence of symmetric commuting extensions for
symmetric matrices may well be provable, but we suspect the question of
minimal dimension is extremely difficult. Fortunately, just because it 
is difficult theoretically does not mean answers cannot be found 
numerically, and we are hopeful that good algorithms can be devised that find 
commuting extensions of a given dimension, if they exist. The determining 
equations are linear and quadratic, and although there surely will be
ill-conditioning in certain cases, it is hard to see why this should be 
so in general. We suspect that the poor performance of the $S(Q,\Lambda_i)$
algorithm in section 3 is to do with the fact we inforced equation
\r{findL}, and only searched over $Q$. Likewise gradient and Newton
algorithms possibly give exaggerated importance to linear components of
the system. There is much more that can be tried here.

Another aspect to be considered
in construction of cubature formulae/commuting extensions  
is symmetry in the domain $\Omega$ and weight $w(x)$. This will clearly 
influence the matrices $A_i$, which represent the natural projections of the
operators $x_i$, and should be respected in construction of the extensions 
$\tilde{A}_i$, see also \cite{cl}. 

We hope very much that more applications will emerge for the notion of
commuting extensions. The idea that noncommutativity can be resolved by
introducing extra dimensions is a very natural one. In fact, we suspect
that, more than the ranks or the norms of commutators, the size of 
minimal commuting extensions is probably the best measure of how 
noncommuting a set of matrices is. 

The minimal size issue appears in 
other settings too. For example, given a set of $m\times n$ matrices $A_i$ we
can ask what is the smallest $N$ such that there exist an $m\times N$
matrix $U$ and an $n\times N$ matrix $V$, both with orthonormal rows, 
and
$N\times N$ diagonal matrices $\Lambda_i$, such that $A_i=U\Lambda_i 
V^T$.
In our context, this provides a natural generalization of the singular
value decomposition of a single matrix, in the same way that (26) 
provides
the generalization of diagonalization of a single symmetric matrix.

\noindent{\bf Acknowledgements:} We thank 
Jonathan Beck, 
Harry Dym, 
Yair Goldfarb, 
Gene Golub, 
Chen Greif, 
Ronny Hadani, 
David Kessler and
Steve Shnider
for suggestions and comments at various stages of this work.

\end{document}